\documentclass{amsart}
\usepackage{amsfonts}
\usepackage{amssymb} 
\usepackage{graphicx}
\usepackage{enumerate}
\usepackage{amsmath}
\usepackage{comment}

\newtheorem{theorem}{Theorem}[section]
\newtheorem{lemma}[theorem]{Lemma}
\newtheorem{corollary}[theorem]{Corollary}

\theoremstyle{definition}

\theoremstyle{remark}

\numberwithin{equation}{section}

\newcommand{\refth}[1]{Theorem~\ref{#1}}
\newcommand{\reflm}[1]{Lemma~\ref{#1}}
\newcommand{\refco}[1]{Corollary~\ref{#1}}

\newcommand{\reffig}[1]{Fig.\ \ref{#1}}
\newcommand{\refeq}[1]{(\ref{#1})}

\newcommand{\beeq}[1]{\begin{equation} \label{#1}}
\newcommand{\eeq}{\end{equation}}
\renewcommand{\(}{\begin{eqnarray*}}
\renewcommand{\)}{\end{eqnarray*}}
\newcommand{\beeqn}{\begin{eqnarray}}
\newcommand{\eeqn}{\end{eqnarray}}
\newcommand{\bs}[1]{{\boldsymbol #1}}

\newcommand{\sothat}{\\ \Rightarrow \hspace*{2mm} &&}
\newcommand{\nexteqline}{\\ &=&}
\newcommand{\eqand}{\mbox{\hspace*{3mm} and \hspace*{3mm}}}
\newcommand{\eqwhere}{\mbox{\hspace*{3mm} where \hspace*{3mm}}}

\renewcommand{\quad}{\hspace*{5mm}}
\renewcommand{\qquad}{\hspace*{10mm}}

\newcommand{\lp}{\left(  }
\newcommand{\rp}{\right) }
\newcommand{\lb}{\left\{  }
\newcommand{\rb}{\right\} }
\newcommand{\lbr}{\left[  }
\newcommand{\rbr}{\right] }
\newcommand{\lf}{\left\lfloor}
\newcommand{\rf}{\right\rfloor}
\newcommand{\lc}{\left\lceil}
\newcommand{\rc}{\right\rceil}

\newcommand{\Z}{{\mathbb Z}}

\newcommand{\R}{{\mathbb R}}

\newcommand{\ep}{\epsilon}
\newcommand{\FF}{{\mathcal F}}

\newcommand{\DD}{{\mathcal D}}
\newcommand{\YY}{{\mathcal Y}}

\newcommand{\TT}{{\mathcal T}}
\renewcommand{\AA}{{\mathcal A}}

\newcommand{\HH}{{\mathcal H}}
\newcommand{\MM}{{\mathcal M}}
\newcommand{\PP}{{\mathcal P}}

\renewcommand{\SS}{{\mathcal S}}

\newcommand{\XX}{{\mathcal X}}

\newcommand{\UU}{{\mathcal U}}

\newcommand{\bfor}{{\bf for }}

\newcommand{\bdo}{{\bf do}:}
\newcommand{\bdoo}{{\bf do }}

\newcommand{\bto}{{\bf to }}

\newcommand{\be}[1]{\begin{enumerate} [#1]}
\newcommand{\ee}{\end{enumerate}}

\newcounter{cnt1}
\newcounter{cnt2}
\newcounter{cnt3}
\newcounter{cnt4}

\newcommand{\bnum}
{
\begin{list}{\arabic{cnt1})}
{
\usecounter{cnt1}
\leftmargin 5mm
\setlength{\leftmargin}{\leftmargin}
\topsep 2pt
\parsep 1pt
\itemsep 1pt}
}
\newcommand{\enum}{\end{list}}

\newcommand{\broman}
{
\begin{list}{\roman{cnt2})}
{
\usecounter{cnt2}
\leftmargin 2mm
\setlength{\leftmargin}{\leftmargin}
\topsep 2pt
\parsep 1pt
\itemsep 1pt}
}
\newcommand{\eroman}{\end{list}}

\newcommand{\bRoman}
{
\begin{list}{\Roman{cnt3})}
{
\usecounter{cnt3}
\leftmargin 5mm
\setlength{\leftmargin}{\leftmargin}
\topsep 2pt
\parsep 1pt
\itemsep 1pt}
}
\newcommand{\eRoman}{\end{list}}

\newcommand{\balph}
{
\begin{list}{\alph{cnt4})}
{
\usecounter{cnt4}
\leftmargin 3mm
\setlength{\leftmargin}{\leftmargin}
\topsep 2pt
\parsep 1pt
\itemsep 1pt}
}
\newcommand{\ealph}{\end{list}}

\newcommand{\bAlph}
{
\begin{list}{\Alph{cnt1})}
{
\usecounter{cnt1}
\leftmargin 5mm
\setlength{\leftmargin}{\leftmargin}
\topsep 2pt
\parsep 1pt
\itemsep 1pt}
}
\newcommand{\eAlph}{\end{list}}

\newcommand{\bdot}
{
\begin{list}{$\cdot$}
{
\leftmargin  3mm
\setlength{\leftmargin}{\leftmargin}
\topsep 3pt
\parsep 1pt
\itemsep 2pt}
}
\newcommand{\edot}{\end{list}}

\newcommand{\bdash}
{
\begin{list}{-}
{
\leftmargin 3mm
\setlength{\leftmargin}{\leftmargin}
\topsep 2pt
\parsep 1pt
\itemsep 1pt}
}
\newcommand{\edash}{\end{list}}

\newcommand{\bnull}
{
\begin{list}{}
{
\leftmargin 3mm
\setlength{\leftmargin}{\leftmargin}
\topsep 0pt
\parsep 1pt
\itemsep 1pt}
}
\newcommand{\enull}{\end{list}}

\begin{document}


\author{Junichiro Fukuyama}
\address{Department of Computer Science and Engineering\\
The Pennsylvania State University\\
PA 16802, USA}
\curraddr{}
\email{jxf140@psu.edu}
\thanks{}




\date{}

\dedicatory{}

\title{Extensions of a Family for Sunflowers}
\begin{abstract} 
This paper explores the structure of the combinatorial domain $2^X$ in relation to sunflowers. The previous study found some intrinsic properties of the {\em $l$-extension}  
\[
Ext \left( \mathcal{F}, l \right) = \left\{ V ~:~ V \in {X \choose l},~ \exists U \in \mathcal{F}~ U \subset V \right\}  
\]
of a family $\FF$ of $m$-cardinality sets. Subsequently, it led to the proof that such an $\mathcal{F}$ includes three mutually disjoint sets if it satisfies the $\Gamma(b)$-condition, that is, 
\(
&& 
\left| \mathcal{F}[S] \right|  < b^{-|S|} |\mathcal{F}|, \quad 
\textrm{for every nonempty set}~ S, 
\\  \textrm{where} && 
\mathcal{F} [S] := \left\{ U : U \in \mathcal{F},~ S \subset U \right\}, 
\)
for $b= m^{\frac{1}{2}+ \epsilon}$ with an $m$ sufficiently larger than a given constant $1/\epsilon$. It is stronger than the statement that $\FF$ includes a 3-sunflower if $|\FF| > b^m$, where {\em $k$-sunflower} refers to a family of $k$ different sets with a common pair-wise intersection. Further refining the theory, we show that an $\FF$ includes $k$ mutually disjoint sets if it satisfies the $\Gamma \left( 8^{\sqrt{\log_2 m}} \sqrt m ~k \log_2 k \right)$-condition with an $m$ sufficiently larger than $k$.  
\end{abstract}
\maketitle

\section{Introduction} In this paper we prove the following theorem\footnote{
An abundance of extra information on the construction is shown in \cite{blog}. Especially, it presents extended proof that such an $\FF$ includes $k$ element-wise disjoint subfamilies $\FF_i$ such that 
\[
|\FF_i| > b^{-\ep m} \lp \frac{\ep}{k \log_2 k} \rp^m |\FF|, \quad \textrm{for every}~ i \in [k], 
\quad \textrm{where}~ \ep = \ln^{-1} c. 
\]
Here two families $\FF_i$ and $\FF_{i'}$ are {\em element-wise disjoint} if $U \cap U' = \emptyset$ for every $U \in \FF_i$ and $U' \in \FF_{i'}$. 
}. 
\begin{theorem} \label{KDC} 
There exists $c \in \R_{>1}$ such that for every $k \in \Z_{>1}$ and $m \in \lp 2^{ck} , \infty\rp \cap \Z$, a family $\FF$ of sets each of cardinality $m$ includes $k$ mutually disjoint sets if it satisfies the $\Gamma(b)$-condition, where 
\[ 
b = c 
2^{\frac{5}{2} \sqrt{\log_2 m}} \lp \sqrt{m} \log_2 m  \rp 
k \log_2 k. 
\qed 
\] 
\end{theorem}

\noindent
Here the {\em $\Gamma(b)$-condition} refers to $|\FF[S]| < b^{-|S|} |\FF|$ for all nonempty sets $S$, where $\FF[S]$ denotes the family of sets in $\FF$ each including the $S$. In particular, the theorem means the existence of $k$ mutually disjoint sets in such an $\FF$ satisfying the $\Gamma\lp b \rp$-condition with $b = 8^{\sqrt{\log_2 m}} \sqrt m ~k \log_2 k$.

\refth{KDC} is stronger than the statement that an $\FF$ with $|\FF| > b^m$ includes a $k$-sunflower, $i.e.$, a family of $k$ different sets with a common pair-wise intersection. As emphasized in \cite{E81}, sunflowers are among the most interested objects in combinatorics since the {\em sunflower lemma} was proven in 1960. The theorem improves the previously known result \cite{3SF} that there exist three mutually disjoint sets in an $\FF$ satisfying the $\Gamma \lp m^{\frac{1}{2} +\ep} \rp$-condition, for any given constant $\ep \in (0, 1)$ and $m$ sufficiently larger than $1/\ep$.

We explore the structure of the combinatorial domain $2^X$ to prove the claim. The previous study revealed some of its intrinsic properties in relation to sunflowers including the following theorem \cite{3SF}. 

\begin{theorem} \label{EGT} 
(Extension Generator Theorem)
There exists $\ep \in (0, 1)$ satisfying the following statement:
let $X$ be the universal set of cardinality $n$, $m \in [n-1]$, $l \in [n]-[m]$,  and
$
\lambda \in \lp 1, \frac{\ep l}{m^2} \rp.
$
For every nonempty family $\FF \subset {X \choose m}$, there exists an $\lp l, \lambda \rp$-extension generator $T$ of $\FF$ with  
$
|T| \le  \lbr \ln {n \choose m} - \ln |\FF| \rbr \Bigr/ \ln \frac{\ep l}{m^2 \lambda}.
$  \qed
\end{theorem}

\noindent 
Here such an {\em $(l, \lambda)$-extension generator $T$} is defined be a set such that 
\[
| Ext \lp \FF[T], l \rp | \ge {n- |T| \choose l-|T|} \lbr 1 - \exp \lp -\lambda \rp \rbr, 
\]
where 
\[
Ext \left( \mathcal{F}, l \right) = \left\{ V ~:~ V \in {X \choose l},~ \exists U \in \mathcal{F}~ U \subset V \right\}, 
\]
is the {\em $l$-extension} of $\FF$. We denote the natural logarithm by $\ln$.

In the process of refining the theory, we confirm \refth{EGT} and other known/new claims in Section 2. We then prove \refth{KDC} in Section 3 by the obtained knowledge.

\section{On the Combinatorial Domain $2^X$} 
This section is organized into four subsections. Sec.\ 2.1 defines basic objects we deal with in this theory to prove a statement called {\em $g^{th}$ mark lemma}. About the norm of a family defined there, Sec.\ 2.2 focuses on a special case proving another lemma. The two claims will be used to prove \refth{KDC}. In Sec.\ 2.3, we define a variation of the $\Gamma$-condition expressed as an inequality, finding some of its structural properties. Lastly, Sec.\ 2.4 develop some methods to partition the universal set $X$ with respect to a given family $\FF$, providing another tool set for our confirmation of \refth{KDC}.

\subsection{Relations between Sets}  
Let $\ep \in (0, 1)$ be sufficiently small depending on no other variables, and $X$ be the universal set of cardinality $n$. Denote 
\(
&& 
\FF, \UU \subset 2^X, \quad 
\\ && 
[b] := \lbr 1, b \rbr \cap \Z, \quad  \textrm{for~} b \in \R_{>0}, 
\\ \textrm{and} && 
{X'\choose m} := \lb U ~:~ U \subset X',~|U| =m  \rb, \quad \textrm{for~} X' \subset X. 
\) 
Also $\FF[S]  = \lb U ~:~ U \in \FF,~ S \subset U  \rb$ as mentioned above, and 
\[
\FF(S) := \lb U ~:~ U \in \FF,~ S \cap U \ne \emptyset  \rb. 
\]

A set means a subset of $X$, and is called {\em $m$-set} if it is in ${X \choose m}$. 
A {\em constant} is a positive real number that depends only on $\ep$. 
A real interval may express the integer set of the same range; for example, $(i, j]$ can denote $(i, j] \cap \Z=\lb i+1, i+2, \ldots, j \rb$. Other basic notations we use in this paper are the same as \cite{3SF}.

Given $g \in \Z_{\ge 2}$, denote $\underbrace{\UU \times \UU \times \cdots \times \UU}_g$ by $\UU^g$ to write 
\[ 
Union \lp \bs{U} \rp := \bigcup_{j=1}^g U_j, \eqand 
Rank \lp \bs{U} \rp := g, 
\quad 
\textrm{for~} \bs{U} = \lp U_1, U_2, \ldots, U_g \rp  \in \UU^g. 
\] 
Let $\FF \subset {X \choose m}$ for some $m \in [n]$, weighted by $w: \lp 2^X \rp^g \rightarrow \R_{\ge 0}$. Define the {\em norm $\left\| \UU \right\|$ of $\UU$} and {\em sparsity $\kappa \lp \FF \rp$} by 
\[
\| \UU \|  = \lbr \sum_{\bs{U} \in \UU^g} w \lp \bs{U}  \rp \rbr^{\frac{1}g}, 
\eqand 
\kappa \lp \FF \rp = \ln {n \choose m}- \ln \| \FF \|, 
\] 
respectively. We may just write $\| U \|$ instead of $\| \lb U \rb \|$ for a single element $U \in \UU$.

For $l \in [n] - [m]$ and $j \in \Z_{\ge 0}$ , let 
\( 
&& 
\PP_{j, g} := \lb \bs{U}~:~ \bs{U} \in \FF^g,~ 
\left| Union\lp \bs{U} \rp  \right| = gm-j 
\rb, 
\\ &&
\DD_g := \lb (\bs{U}, Y) ~:~ \bs{U} \in \FF^g,~ 
Y \in {X \choose l},~ Union(\bs{U}) \subset Y 
\rb, 
\\ &&
\| \PP_{j, g} \| := \sum_{\bs{U}\in \PP_{j, g}} w  \lp \bs{U} \rp,
\eqand 
\| \DD_g \| := \sum_{(\bs{U}, Y) \in \DD_g} w  \lp \bs{U} \rp, 
\)
extending the norm $\| \cdot \|$ for $\PP_{j, g}$ and $\DD_g$. The weight function $w$ is said to {\em induce} the norm $\| \cdot \|$ and sparsity $\kappa$.

The family $\FF$ satisfies the {\em $\Gamma_g (b, h)$-condition on $\| \cdot \|$} $\lp b, h  \in \R_{>0} \rp$  if 
\beeqn 
&& \label{eqGamma_g}
\| \UU \|  = \lbr \sum_{\bs{U} \in \lp \UU \cap \FF \rp^g} w \lp \bs{U}  \rp \rbr^{\frac{1}g}, \quad 
\textrm{for all~} \UU \subset 2^X, 
\\ \textrm{and} && \nonumber 
\| \PP_{j, g} \| 
< h b^{-j} 
\| \FF \|^g, \quad 
\textrm{for every~} j \in [(g-1)m]. 
\eeqn 
We may drop the subscript $g$ if it is clear from the context, so $\PP_{j, g}$ can be written as $\PP_j$, $\Gamma_g (b,  h)$ as $\Gamma (b, h)$ etc. Here are the statement mentioned above that is proven in \cite{3SF} and Appendix A.

\begin{lemma} \label{gthMark} 
Let 
\be{i)} 
\item $X$ be the universal set weighted by $w: \lp 2^X \rp^g \rightarrow \R_{\ge 0}$ for some $g \in \Z_{\ge 2}$, 
\item $l \in [n]$, $m \in [l-1]$, and $h, \gamma \in \R_{>0}$, such that 
$\gamma$ and $\frac{l}{gm}$ are both sufficiently large, 
\item and $\FF \subset {X \choose m}$ satisfy the $\Gamma_g  \lp \frac{4 \gamma n}{l},~h \rp$-condition on the norm $\| \cdot \|$ induced by $w$.  
\ee 
Then 
\[ 
\| \DD_g  \| < 
\frac{\lp 1 + \frac{h}{\gamma} \rp {n \choose l}{l \choose m}^g}{{n \choose m}^g} 
\| \FF \|^g 
.
\qed 
\]
\end{lemma} 

\medskip 

Call the statement {\em $g^{th}$ mark lemma} as an element $(\bs{U}, Y)$ of $\DD_g$ is said to be a {\em $g^{th}$ mark} of $\FF$ as in \cite{blog}. 

Since 
\[
\| \DD_g \| 
= \sum_{(\bs{U}, Y) \in \DD_g} w \lp \bs{U} \rp
=\sum_{Y \in{X \choose l}} 
~\sum_{\bs{U} \in \lbr \FF \cap {Y \choose m} \rbr^g} 
w \lp \bs{U} \rp 
=\sum_{Y \in{X \choose l}}  \left\|  {Y \choose m} \right\|, 
\]
it means: 

\begin{corollary} \label{EGTTildeG}
For such objects and $\ep \in (0, 1)$, there are $\lc \lp 1 - \ep \rp {n \choose l} \rc$ sets $Y \in {X \choose l}$ such that
\[
\left\| {Y \choose m} \right\| 
< 
\frac{\lp 1 + \frac{h}{\gamma} \rp {l \choose m}^g}{\ep {n \choose m}^g}
\left\| \FF \right\|. 
\qed
\]
\end{corollary}

\subsection{Primitive Weight on $X$} 
Now let the weight $w$ be defined by 
\[
(U_1, U_2, \ldots, U_g) \mapsto \prod_{i=1}^g w_*(U_i), 
\quad 
\textrm{for some~}
w_*: 2^X \rightarrow \R_{\ge 0}. 
\] 
We say that such an $X$ is {\em primitively weighted} by $w$ and $w_*$, and that the weight $w$, the induced norm $\| \cdot \|$ and sparsity $\kappa$ are {\em primitive with $w_*$}.

We have the following theorem.

\begin{theorem} \label{EGT42} 
Let $X$ be primitively weighted to induce the norm $\| \cdot \|$. 
For every sufficiently small $\ep \in (0, 1)$, and $\FF \subset {X \choose m}$ satisfying the $\Gamma_2 \lp \frac{4 \gamma n}{l},~1  \rp$-condition on $\| \cdot \|$ for some $l \in [n]$, $m \in [l]$, and 
$
\gamma \in \lbr \ep^{-2},~ l m^{-1} \rbr 
$, there are $\lc {n \choose l} \lp 1-  \ep \rp \rc$ sets $Y \in {X \choose l}$ such that
\[
\lp 1 - \sqrt{\frac{2}{\ep \gamma}} \rp
\frac{{l \choose m}}{{n \choose m}}
\left\| \FF \right\|
<
\left\| {Y \choose m} \right\|
<
\lp 1 + \sqrt{\frac{2}{ \ep \gamma}}  \rp
\frac{{l \choose m}}{{n \choose m}}
\left\| \FF \right\|.
\qed
\]
\end{theorem}

\medskip 

\noindent 
It is proven in Appendix A as well as \cite{blog, 3SF} by the $g^{th}$ mark lemma for $g=2$. The theorem can be used to prove the extension generator theorem as demonstrated in Appendix A.

Let us say that an $\FF$ satisfies the {\em $\Gamma(b', h')$-condition on $\| \cdot \|$}  $(b', h' \in \R_{>0})$, if 
$\|\FF [S] \| < h'  b'^{-|S|} \| \FF \|$ for any nonempty set $S$. We may constrain the range of $S$ to be $2^{X'}$ for $X' \subset X$ to say the {\em $\Gamma(b', h')$-condition on $\| \cdot \|$ in $X'$}. In addition, the parameter $h'$ could be skipped if it is 1.

Remarks.

\be{A)} 
\item For a primitive norm, 
\[
\left\| \UU \right\|  
= \sum_{U \in \UU} 
w_* \lp U \rp, 
\quad 
\textrm{for any~} \UU \subset 2^X, 
\] 
since 
$
\lbr \sum_{U \in \UU 
} 
w_* \lp U \rp \rbr^g 
= 
\sum_{\bs{U} \in \UU^g} w \lp \bs{U} \rp
= 
\| \UU \|^g 
$. By this {\em linearity of primitive norm $\| \cdot \|$}, the first condition of \refeq{eqGamma_g} for $\FF$ on $\| \cdot \|$ is equivalenet to $\| \UU \| = \| \UU \cap \FF \|$ for all $\UU \subset 2^{X}$.

\item Suppose an $\FF$ satisfies the $\Gamma(b m)$-condition on a primitive norm $\| \cdot \|$. If, in addition, the first condition of \refeq{eqGamma_g} holds for $\FF$, it satisfies the $\Gamma_2(b, 1)$-condition. It is true since for each $u \in [m]$, 
\(
\sum_{U_1, U_2 \in \FF \atop |U_1 \cap U_2| = u} w(U_1, U_2) &=& 
\sum_{U_1, U_2 \in \FF \atop |U_1 \cap U_2| = u} w_*(U_1) w_*(U_2) 
\\ &\le&  
\sum_{U_1 \in \FF \atop U \in {U_1 \choose u}} w_*(U_1) \| \FF [U] \| 
\\ &\le&  
\| \FF \|^2 (bm)^{-u} {m \choose u}
\le b^{-u} \| \FF \|^2. 
\)
The second line holds by the first condition of \refeq{eqGamma_g}, and the last line by the $\Gamma(bm)$-condition of $\FF$ on $\| \cdot \|$. 

\item The $\Gamma(b)$-condition on $\| \cdot \|$ induced by the unit weight 
\[
w_* (U) = 
\lb \begin{array}{cc}
1, &  \textrm{if $U \in \FF$,} \\
0, &  \textrm{otherwise,} \\
\end{array} \right. 
\] 
is identified with the $\Gamma(b)$-condition referred to by \refth{KDC}.

\item We have $\kappa \lbr {X \choose 2m} - Ext \lp \FF, 2m \rp \rbr \ge 2 \kappa \lbr {X \choose m} - \FF \rbr$ for the sparsity $\kappa$ induced by the unit weight. It is shown in \cite{blog, 3SF} and by \reflm{PhaseII} in Appendix A.

\item In addition, $\kappa \lbr Ext \lp \FF, l \rp \rbr \le \kappa \lp \FF \rp$ for $l \in [n]-[m]$. For there are $|\FF| {n-m \choose l-m} = {n \choose m}e^{-\kappa \lp \FF \rp} {n-m \choose l-m}= {n \choose l}{l \choose m} e^{-\kappa \lp \FF \rp}$ set pairs $(S, T)$ such that $S \in \FF$, $T \in {X \choose l}$ and $S \subset T$. This means $|Ext \lp \FF, l \rp| \ge {n \choose l}e^{-\kappa \lp \FF \rp}$ leading to the claim. \qed 
\ee 

\medskip

By \refth{EGT42} and B): 

\begin{corollary} \label{corEGT42} 
Let $\ep$, $l$, $m$ and $\gamma$ be as given by \refth{EGT42}. For a family 
$\FF \subset {X \choose m}$ satisfying the $\Gamma \lp \frac{4 \gamma n m}{l} \rp$-condition on a primitive norm $\| \cdot \|$, there exist $\lc {n \choose l} \lp 1-  \ep \rp \rc$ sets $Y \in {X \choose l}$ such that
\[
\lp 1 - \sqrt{\frac{2}{\ep \gamma}} \rp
\frac{{l \choose m}}{{n \choose m}}
\left\| \FF \right\|
<
\left\| \FF \cap {Y \choose m}  \right\|
<
\lp 1 + \sqrt{\frac{2}{ \ep \gamma}}  \rp
\frac{{l \choose m}}{{n \choose m}}
\left\| \FF \right\|.
\qed
\]
\end{corollary}

The following lemma also verified in Appendix A will be used by our proof of \refth{KDC}.

\begin{lemma} \label{lmConversion} 
(Conversion Lemma) 
Let $X$ be weighted by $w: \lp 2^X \rp^g \rightarrow \R_{\ge 0}$ inducing $\| \cdot \|$, primitively with some norm $\| \cdot \|_*$. If a family $\FF \subset {X \choose m}$ satisfies the $\Gamma(b, h)$-condition on $\| \cdot \|_*$ for some $b, h \in \R_{\ge 1}$, it satisfies the $\Gamma_g \lbr \frac{b}{2^{g-2} (g-1) m},~ h^{g-1} \rbr$-condition on $\| \cdot \|$. 
\qed 
\end{lemma}

\subsection{Collective $\Gamma$-Conditions}
With the two lemmas \ref{gthMark} and \ref{lmConversion}, we study the {\em $g^{th}$ collective $\Gamma(b, h)$-condition of $\FF \subset {X \choose m}$ on $\| \cdot \|$} in this subsection. It is defined to be 
\beeq{eqCGamma}
\sum_{u \in [m] \atop U \in {X \choose u}} \frac{\| \FF[U] \|^g}{b^{-u(g-1)} {m \choose u}}
< h \| \FF \|^g, 
\eeq 
for $g \in \Z_{\ge 2}$, $b, h \in \R_{>0}$ and a primitive norm $\| \cdot \|$.

It is straightforward to check that: 
\be{A)} \setcounter{enumi}{5}
\item If $\FF$ satisfies the $\Gamma(2 b)$-condition on $\| \cdot \|$, it satisfies the $g^{th}$ collective $\Gamma(b, 1)$-condition. 
\item If $\FF$ satisfies the $g^{th}$ collective $\Gamma(b, 1)$-condition on $\| \cdot \|$, it satisfies the $\Gamma(b m^{-3/g})$-condition. 
\qed 
\ee 
They are shown in \cite{blog} as well. In addition, as confirmed in Appendix A:

\begin{lemma} \label{pushup}
(Push-Up Lemma) 
In the universal set $X$ weighted with a primitive norm $\| \cdot \|$, let $\FF \subset {X \choose m}$ satisfy 
\[
\sum_{u \in [u_*, m] \atop S \in {Z \choose u}
}
\frac{\|\FF[S]\|^g}{b^{-(g-1)u} {m \choose u}}
< h^2 \|\FF \|^g, 
\]
for some set $Z \subset X$, integers $g \in \Z_{\ge 2}$ and $u_* \in [0, m]$, and real numbers $b \in \R_{>1}$ and $h \in (0, 1/4]$. Then there exists $S \in {Z \choose v}$ for some $v \in [0, u_*)$ satisfying the three conditions. 
\( 
1) && 
\sum_{u \in [u_* -v, m-v] \atop 
T \in {Z-S \choose u} 
}
\frac{\|\FF[S \cup T]\|^g}{b^{-(g-1)u} {m-v \choose u}}
< h \|\FF[S] \|^g. 
\\ 2) && 
\sum_{u \in [u_* - v] \atop 
T \in {Z-S \choose u} 
}
\frac{\|\FF[S \cup T] \|^g}{\lp \frac{b}{6} \rp^{-(g-1)u} {m-v \choose u}}
< \|\FF[S] \|^g. 
\\ 3) && 
\|\FF[S] \| >  h^{\frac{1}{g-1}} b^{-v} \|\FF\|. \qed
\)
\end{lemma}


\subsection{Partitioning $X$ by a Given $\FF$} 
In this subsection, we develop tools to partition the universal set $X$ with respect to a particular family of $m$-sets. The three theorems presented below are all proven in Appendix B.  Given an $r \in [n]$, a {\em split of $X$} is a tuple $\bs{X}=(X_1, X_2, \ldots, X_r)$ of $r$ mutually disjoint nonempty sets $X_i$, called {\em strips} of $\bs{X}$, whose union is $X$.  We have a general fact here.

\medskip

\begin{theorem} \label{thSplit} 
Let 
\bdash 
\item a universal set $X$ of cardinality $n$ be weighted with a primitive norm $\| \cdot \|$, 
\item $m \in [n]$, and $r \in [m]$, 
\item $d_i \in [n]$ for $i \in [r]$ with $\sum_{i=1}^r d_i= n$, and $q_i \in [d_i]$ with $\sum_{i=1}^r q_i= m$. 
\edash 
For a family $\FF \subset {X \choose m}$, there exists a split $\bs{X}=(X_1, X_2, \ldots, X_r)$ of $X$ and subfamily $\FF' \subset \FF$ such that  
\[
\|\FF' \| \ge \frac{\|\FF \|}{{n \choose m}} \prod_{i=1}^r {d_i \choose q_i}, 
\quad 
|X_i|= d_i, 
\eqand 
|U \cap X_i| = q_i, 
\]
for all $i \in [r]$ and $U \in \FF'$. \qed 
\end{theorem}

\medskip

Let $r$ divide $n$. An {\em $r$-split} is a split $\bs{X}$ such that $|X_i| = n / r$ for all $i \in [r]$. An $m$-set $U$ is said to be {\em on the $\bs{X}$} if $|X_i \cap U| = \lf m/r \rf$ for $i <r$. Denote by ${\bs{X} \choose m}$ the family of $m$-sets on $\bs{X}$. The theorem implies:

\begin{corollary} \label{Split2} 
For a family $\FF \subset {X \choose m}$ in such an $X$, and $r \in [m]$ that divides both $m$ and $n$, there exists an $r$-split $\bs{X}$ of $X$ such that 
\[
\left \| \FF \cap {\bs{X} \choose m} \right \| \ge \frac{{n/r \choose m/r}^r}{{n \choose m}}  \| \FF \|. \qed 
\]
\end{corollary}

As $\lp \frac{n}{m} \rp^m > {n \choose m} \exp \lp -m \rp$ by \refeq{eq1App} in Appendix C:

\begin{corollary} \label{Split3} 
In a universal set $X$ primitively weighted with an induced sparsity $\kappa$, any family $\FF \subset {X \choose m}$ with $m \big| n$ and finite $\kappa \lp \FF \rp$, there exists an $m$-split $\bs{X}$ of $X$ such that 
$
\kappa \lp \FF \cap {\bs{X} \choose m} \rp  
<
\kappa \lp \FF \rp + m
$. \qed  
\end{corollary}

If we relax the condition $|X_i \cap U| = \frac{m}{r}$ to $\frac{m}{2r}< |X_i \cap U| < \frac{2m}{r}$ for $U \in \FF'$, the sparsity increase $\kappa \lp \FF' \rp- \kappa \lp \FF \rp$ can be seen much smaller than a constant. The following statement connfirms it. 

\begin{theorem} \label{thSplit2} 
Let $X$ be weighted to induce a primitive norm $\| \cdot \|$, 
$m \in [n]$, and $r \in [m]$ such that $r \big| n$ and $\max \lp \frac{m}{n}, \frac{r}{m}, \frac{1}{r} \rp< \ep$. For a family $\FF \subset {X \choose m}$ with $\| \FF \| >0$, there exists an $r$-split $\bs{X}=(X_1, X_2, \ldots, X_r)$ and subfamily $\FF' \subset \FF$ satisfying 
\[
\|\FF' \| > \lbr 1-  r m \exp \lp - \frac{m}{12r} \rp \rbr \| \FF \|, \eqand \frac{m}{2r} < |U \cap X_i| < \frac{2m}{r}, 
\]
for all $U \in \FF'$ and strips $X_i$ of $\bs{X}$. \qed 
\end{theorem}

\medskip

By the theorem, we can show that:

\begin{theorem} \label{thSplit3} 
Let $m \in [n]$, $r \in [m]$ such that 
\[ 
r \big| n, 
\quad 
\max \lp \frac{m}{n}, \frac{r}{m}, \frac{1}{r} \rp< \ep, \eqand b \in  \lp c^3 kr \ln \frac{m}{r}, \infty \rp, 
\]
for some constants $c \in ( \ep^{-1}, \infty)$ and $k \in \Z_{>0}$. For each nonempty $\FF \subset {X \choose m}$ satisfying the $\Gamma(b)$-condition, there exist an $r$-split $\bs{X}=(X_1, X_2, \ldots, X_r)$, $k$ mutually disjoint sets $C_i$, and subfamilies $\FF_i \subset \FF[C_i]  \cap {X - \bigcup_{i' \in [k] - \lb i \rb} C_{i'} \choose m}$ with the following four conditions satisfied for every $i \in [k]$: 
\be{i)}
\item $|C_i| < c r^2 \ln \frac{m}{r}$. 
\vspace*{1.2mm}
\item $|\FF_i|> \frac{1}{2} \lp \frac{r}{3m} \rp^r  b^{-|C_i|} |\FF|$. 
\item The $\Gamma\lp b / 6 \rp$-condition of $\FF_i$ in $X - C_i$. 
\item $|U \cap X_j| = q_{i, j}$ for every strip $X_j$ of $\bs{X}$, some $q_{i, j} \in \lp \frac{m}{2r},~ \frac{2m}{r}  \rp$, and all $U \in \FF_i$. \qed 
\ee
\end{theorem}

Our proof of \refth{KDC} in the next section will refer to \refth{thSplit3}.

\section{Proof of \refth{KDC}} 
Let a constant $c$ be larger than $2^{1 / \ep}$. 
Given such a $k$ and $\FF \subset {X \choose m}$ by the statement, let 
\[
b = 2^c k_\dag m_\dag, \eqwhere  k_\dag = k \lg k,  \eqand m_\dag =  2^{\frac{5}{2} \sqrt{\lg m}} \sqrt{m} \lg m, 
\] 
where $\lg$ denotes $\log_2$. 
WLOG assuming the $\Gamma(b)$-condition of $\FF$ with this $b$, we will find $k$ mutually disjoint sets in $\FF$ as our proof of \refth{KDC}.

Given an $m> 2^{c k}$ as well, put 
\[ 
g = \lf \sqrt{\lg m} \rf,  \eqand 
r = \lf \frac{1}{c^2}  \sqrt{\frac{m}{\lg m}} \rf. 
\]
We also assume $n=|X|$ is sufficiently larger than $m$ divisible by $mr$; otherwise, add some extra elements to $X$ to make this true.

By \refth{thSplit3}, there exist an $r$-split $\bs{X} = (X_1, X_2, \ldots, X_r)$, mutually disjoint sets $C_i$, and subfamilies $\FF_i \subset \FF[C_i] - {X - \bigcup_{i' \in [k] - \lb i \rb} C_{i'} \choose m}$ $\lp i \in [k] \rp$ such that i)-iv). In paricular, we have the integers $q_{i, j} \in \lp \frac{m}{2r}, \frac{2m}{r} \rp$ meeting iv) that appear explicitly in our construction below.

Update $b$ by $b \leftarrow b/12$, so we have the $g^{th}$ collective $\Gamma$-condition 
\beeq{eq1Main}
\sum_{u \in [m] \atop S \in {X-C_i \choose u}} \frac{|\FF_i[S]|^g}{b^{-(g-1)u} {m \choose u}} 
< |\FF_i|^g. 
\eeq 
Its truth is due to Remark F) of Sec.\ 2.3.

\medskip

With the obtained objects we consider the following statement for all $j \in [0, r-1]$: 

\medskip 

\noindent 
{\bf Proposition $\Pi(j)$:} there exist $k$ mutually disjoint sets $C_i$ and subfamilies $\FF_i \subset \FF[C_i] \cap {X - \bigcup_{i' \in [k] - \lb i \rb} C_{i'} \choose m}$ satisfying the four conditions for every $i \in [k]$: 
\(
\textrm{i)} &&  
|C_i| < \frac{m}{c} + j r_*, \eqwhere r_* := \frac{m}{c r}. 
\\ \textrm{ii)} && 
\sum_{u \in [r_*, m_i] \atop T \in {Z_{i, j} \choose u}}~ \frac{|\FF_i[T]|^g}{b_j^{-(g-1)u} {m_i \choose u}} 
< |\FF_i|^g, 
\\ && 
\textrm{where}  \qquad 
m_i := m - |C_i|, \quad 
Z_{i, j} := X- C_i - \bigcup_{j'=1}^j X_{j'}, 
\\ && 
\qquad \qquad \textrm{and}  \qquad 
b_j := b \lp 1- \frac{1}{r} \rp^{j}, 
\\ \textrm{iii)} && 
\sum_{u \in [m_j] \atop T \in {Z_{i, j} \choose u}}~ \frac{|\FF_i[T]|^g}{\lp \frac{b_j}{6} \rp^{-(g-1)u} {m_i \choose u}} 
< 2 |\FF_i|^g. 
\\ \textrm{iv)} && 
\textrm{$\FF_i$ is {\em element-wise disjoint} with the other $\FF_{i'}$ within $\bigcup_{j'=1}^j X_{j'}$, $i.e.$,} 
\\ && 
\quad 
U \cap U' \cap \bigcup_{j'=1}^j X_{j'}= \emptyset, 
\quad \forall U \in \FF_i  \hspace*{2mm} \forall U' \in \bigcup_{i' \in [k] - \lb i \rb} \FF_{i'}. 
\qed 
\)


\medskip 

We prove it by induction on $j$. After its completion, we will find that $\Pi(r-1)$ implies the existence of $k$ mutually disjoint sets in $\FF$. To start our proof, we verify $\Pi(0)$ as the induction basis: choose $C_i$ and $\FF_i$ to be the ones above from \refth{thSplit3}. Then $\Pi(0)$-i) and iv) are clear in addition to $\Pi(j)$-ii) and iii) from \refeq{eq1Main}.

Assume $\Pi(j)$ for $j \le r-2$. When we consider each $C_i$, write for simplicity 
\(
&& 
m_* := q_{i, j+1} - |C_i \cap X_{j+1}|, \quad 
X_* := X_{j+1} - \bigcup_{i \in [k]} C_i, 
\\ && 
\HH := {X_* \choose m_*},  \quad 
n_* := |X_*|,   \quad  
l  := \lf \frac{\ep n_*}{k_\dag} \rf. 
\) 
Assume $m_*>0$ for all $i \in [k]$ for a while. We prove $\Pi(j+1)$ in five steps. 

\medskip

{\bf Step 1.~}{\em Construct a family $\YY_i$ of $Y \in {X_* \choose l}$ such that $\FF_i \cap {X - X_* \cup Y \choose m}$ is sufficiently \\[1ex] large. } Consider each $i \in [k]$. We claim that $\FF_i$ satisfies the $\Gamma\lp c^6 g^6 k_\dag m_* \rp$-condition in $Z_{i, j}$. For, by $\Pi(j)$-iii), 
\[
|\FF_i[T]| < \lbr \lp \frac{b_j}{12} \rp^{1- \frac{1}{g}} m^{-\frac{1}{g}} \rbr^{-u}  |\FF_i|, \quad 
\forall u \in [m_i] ~ 
\forall T \in {Z_{i, j} \choose u}. 
\]
Here 
\(
&& 
b_j \ge  b \lp 1- r^{-1} \rp^{r} > b \exp \lp - 2 \rp> 2^{c/2} k_\dag m_\dag, 
\sothat 
\lp \frac{b_j}{12} \rp^{1- \frac{1}{g}} m^{-\frac{1}{g}} > 
2^{c/3}  k_\dag m_\dag  (k_\dag m_\dag m)^{- \frac{1}{g}} > 2^{c/4} g^7 k_\dag \sqrt m, 
\)
where the last inequality is due to $m_\dag < 8^{\sqrt{\lg m}} \sqrt m$ and $m > 2^{c k}$ implying 
\(
\lp k_\dag m_\dag m \rp^{- \frac{1}{g}} &>& \lp k^2 8^{\sqrt{\lg m}}  m^{\frac{3}{2}} \rp^{-1 \big/ \lf \sqrt{\lg m} \rf} 
> 2^{-(1+\ep)\frac{3}{2} \sqrt{\lg m}}
\\ &>& g^7 2^{- \frac{5}{2} \sqrt{\lg m}}. 
\) 
By those and 
\beeq{eq2Step2}
m_* \le q_{i, j} < \frac{2m}{r} < c^3 \sqrt{m \lg m}, 
\eeq
we have the claimed $\Gamma\lp c^6 g^6 k_\dag m_* \rp$-condition of $\FF_i$. Note it also satisfies the \\$\Gamma(c^6 g^6 k_\dag \sqrt m)$-condition we will refer to later.

Weight $X_*$ by $w : \lp 2^{X_*} \rp^2  \rightarrow \Z_{\ge 0}$ primitively with 
$
U \mapsto |\FF [U]|
$   
inducing the norm $\| \cdot \|$. Then $\HH$ satisfies the $\Gamma(c^6 g^6 k_\dag m_*)$-condition on $\| \cdot \|$ by the claim. 

Apply \refco{corEGT42} to $\HH$ with the detected properties. There exists a family $\YY_i \subset {X_* \choose l}$ such that
$
|\YY_i| > {n_* \choose l} \lp 1  - \frac{1}{cg} \rp
$,
and 
\beeqn 
&& \label{eq2Step1}
\frac{{l \choose m_*}}{{n_* \choose m_*}} |\FF|
\lp 1  - \frac{1}{cg} \rp
< 
\left| \FF_{i, Y} \right|
< 
\frac{{l \choose m_*}}{{n_* \choose m_*}} |\FF|
\lp 1 + \frac{1}{cg} \rp, \quad 
\forall Y \in \YY_i, 
\\ && \nonumber 
\textrm{where} \quad \FF_{i, Y} = \FF \cap {X - X_* \cup Y \choose  m}.
\eeqn

\medskip

{\bf Step 2.~}{\em Show that most $Y \in \YY_i$ satisfies $\Pi(j+1)$-ii) with the substitution $\FF_i \leftarrow \FF_{i, Y}$.} 
Keep fixing each $i$ to write 
\( 
&& 
\SS = \bigcup_{u=0}^{m_*} {X_* \choose u}, \quad 
S \in \SS, \quad 
Z = X - C_i - \bigcup_{j'=1}^{j+1} X_{j'}, 
\\ && 
\TT = \bigcup_{u=[r_*, m_i]} {Z \choose u}, \quad 
T \in \TT. 
\) 
Abusing the symbols $\FF$, $m$ and $b$, we may also denote 
\[
\FF = \FF_i, \quad m=m_i, \quad b = b_j, 
\]
for simplicity if clear from context.

Define 
\( 
w_T: \lp 2^{X_*} \rp^g \rightarrow \R_{\ge 0}, \quad 
(U_1, U_2, \ldots, U_g) \mapsto 
\frac{\prod_{i=1}^g \left| \FF[U_i \cup T] \right|}
{
b^{-(g-1)|T|} {m \choose |T|} 
}, 
\) 
for each $T \in \TT$ inducing the norm $\| \cdot \|_T$. Reset the weight $w$ to  
\(
w: \lp 2^{X_*} \rp^g \rightarrow \R_{\ge 0}, \quad 
\bs{U} \mapsto
\sum_{T \in \TT}  
w_T \lp \bs{U} \rp, 
\)  
inducing the norm $\| \cdot \|$. For each $T$, 
\[
w_{*, T}: 2^{X_*}  \rightarrow \R_{\ge 0}, \quad 
U  \mapsto  
\frac{\left| \FF[U \cup T] \right|}
{
b^{-\lp 1- \frac{1}{g} \rp|T|} {m \choose |T|}^{\frac{1}{g}} 
}, 
\] 
inducing the norm $\| \cdot \|_{*, T}$, and 
\[
h_T = 
\| \HH \|_{*, T}^{-1}
\lbr 
\sum_{S \in \SS}  
\frac{\left| \FF[S \cup T] \right|^g}
{
b^{-(g-1) \lp s + |T| \rp} {m  \choose {s + |T|}}
}
\rbr^{\frac{1}{g}}. 
\]
We perform the task in two substeps. 

\medskip 

{\em Step 2-1.~}{\em Confirm some $\Gamma_g$-condition of $\HH$.} 
Below we show the $\Gamma_g \lp b_\circ, h \rp$ condition of $\HH$ on $\| \cdot \|$, where 
\[
b_\circ = \frac{b_\dag}{2^g (g-1) m_*} , \eqand h = \frac{|\FF|^g}{\| \HH \|^g}, 
\eqwhere
b_\dag = \frac{b^{1- \frac{1}{g}}}{m^{\frac{1}{g}}}. 
\]

\medskip 

Remarks. 

\be{A)} 
\item 
\[
\sum_{S \in \SS \atop T \in \TT}~ \frac{|\FF[S \cup T]|^g}{b^{-(g-1)(|S|+|T|)} {m \choose |S|+|T|}} 
< |\FF|^g, 
\]
holds due to $\Pi(j)$-ii).

\item  \label{rem1} For each $T \in \TT$, the norm $\| \cdot \|_T$ is primitive with $\| \cdot \|_{*, T}$ satisfying 
\[
\| \HH[S] \|_{*, T}^g = 
\| \HH[S] \|_T^g = 
\sum_{(U_1, U_2, \ldots, U_g) \in \HH[S]^g} 
\frac{\prod_{i=1}^g \left| \FF[U_i \cup T] \right|}
{b^{-(g-1)|T|} {m \choose |T|}}
= \frac{\left| \FF[S \cup T] \right|^g}
{b^{-(g-1)|T|} {m \choose |T|}}, 
\]
for all $S \in \SS$.

\item By B), 
\[
\sum_{S \in \SS}  
\frac{\| \HH[S] \|_{*, T}^g}{b^{-(g-1)s}~{m \choose s} } 
\le 
\sum_{S \in \SS}  
\frac{\left| \FF[S \cup T] \right|^g}
{
b^{-(g-1)\lp s + |T| \rp} {m  \choose s + |T|}
}
= 
h_T^g \| \HH \|^g_{*, T}, 
\]
since 
\[
{m \choose s}  {m \choose |T|} 
\ge 
{m \choose s} {m -s \choose s+|T| - s} 
= 
{m \choose s+|T|}{s+|T| \choose s}
\ge 
{m \choose s+|T|}. 
\]
The leftmost inequality is strict if $s$ is positive. 

\item Thus, 
\[ 
\| \HH[S] \|_{*, T} < 
h_T \| \HH \|_{*, T}  b^{ - \lp 1- \frac{1}{g} \rp |S|} m^\frac{|S|}{g}
= h_T b_\dag^{-|S|} \| \HH \|_{*, T}. 
\]
for each $T$ and nonempty $S$, that is, $\HH$ satisfies the $\Gamma(b_\dag, h_T)$-condition on $\| \cdot \|_{*, T}$. 

\item $h_T \ge 1$ since 
\( 
&& 
h_T^g = 
\| \HH \|_{*, T}^{-g} 
\sum_{S \in \SS}  
\frac{\left| \FF[S \cup T] \right|^g}
{
b^{-(g-1) \lp s + |T| \rp} {m  \choose {s + |T|}}
} \ge 1, 
\\ \textrm{from} &&
\| \HH \|_{*, T} = \sum_{U \in \HH} \frac{\left| \FF[U \cup T] \right|}
{
b^{-\lp 1- \frac{1}{g} \rp|T|} {m \choose |T|}^{\frac{1}{g}} 
}
= \frac{\left| \FF[T] \right|}
{
b^{-\lp 1- \frac{1}{g} \rp|T|} {m \choose |T|}^{\frac{1}{g}} 
}. 
\) 
By the conversion lemma in Sec.\ 2.2, $\HH$ satisfies the $\Gamma_g (b_\circ, h_T^g)$ condition on $\| \cdot \|_T$ for each $T$. 
\qed 
\ee 

\medskip 

We show the desired $\Gamma_g$-condition of $\HH$ on $\| \cdot \|$ by the remarks. By the derived $\Gamma_g$-condition of $\HH$ for each $T$, 
\[
\| \PP_{j', g}  \|_T <  h_T^g b_\circ^{-j'} \| \HH \|_T^g  = 
b_\circ^{-j'} \sum_{S \in \SS}  
\frac{\left| \FF[S \cup T] \right|^g}
{
b^{-(g-1) \lp s + |T| \rp} {m  \choose {s + |T|}}
}, 
\]
for every $j' \in [(g-1) m_*]$. 
Sum this up for all $T \in \TT$ to see 
\[
\| \PP_{j', g}  \| < b_\circ^{-j'} \sum_{S \in \SS \atop T \in \TT}  
\frac{\left| \FF[S \cup T] \right|^g}
{
b^{-(g-1) \lp s + |T| \rp} {m  \choose {s + |T|}}
}
< b_\circ^{-j'} |\FF|^g 
= h b_\circ^{-j'} \| \HH \|^g, 
\] 
by A). 
Hence, $\HH$ satisfies the $\Gamma_g(b_\circ, h)$-condition on $\| \cdot \|$.

\medskip

{\em Step 2-2.~}{\em Eliminate $Y \in \YY_i$ with excessively large 
$\sum_{T \in \TT} |\FF_{i, Y}[T]|^g$.} 
Regarding the confirmed $\Gamma_g$-condition, we note that 
\(
&& 
2^g  (b m)^\frac{1}{g} < 2^g \lp 8^{\sqrt{\lg m}} m^{\frac{3}{2}}\rp^\frac{1}{\lf \sqrt{\lg m} \rf}
<2^{4+\frac{5}{2} \sqrt{\lg m}}, 
\eqand \textrm{\refeq{eq2Step2} in Step 1}, 
\sothat 
b_\circ = \frac{b^{1- \frac{1}{g}}}{g 2^g m_* m^{\frac{1}{g}}}> c^4 k_\dag.   
\)
In addition, $h> 1$; otherwise,  
$
|\FF|^g \le \| \HH  \|^g 
=
\sum_{T \in \TT}  
\frac{|\FF[T]|^g  
}
{
b^{-(g-1) |T|} {m_* \choose |T|}
} 
< 
|\FF|^g 
$ 
by A).

Apply \refco{EGTTildeG} to $\HH$ on $\| \cdot \|$. There are more than $\lp 1-   c^{-1} \rp {n_* \choose l}$ sets $Y \in {X_* \choose l}$ such that
\(
&&
\sum_{T \in \TT } ~
\frac{|\FF_{i, Y}[T]|^g  
}
{
b^{-(g-1)|T|} {m \choose |T|}
}
= 
\sum_{\bs{U} \in \lbr \HH \cap {Y \choose m_*} \rbr^g} w \lp \bs{U} \rp 
\\ &<&
\frac{c^2 \lp 1 + \ep h \rp {l \choose m_*}^g}{{n_* \choose m_*}^g}
\| \HH \|^g 
< 
\frac{c^2 (1+\ep) {l \choose m_*}^g}{{n_* \choose m_*}^g}
|\FF|^g
\\ &<& \nonumber 
2 c^2 |\FF_{i, Y}|^g, 
\) 
where the last inequality is due to \refeq{eq2Step1}. It means 
\beeq{eqStep2_2}
\sum_{
u \in [r_*, m_i] 
\atop 
T \in {Z \choose u}
}~ 
\frac{|\FF_{i, Y}[T]|^g  
}
{
b_j^{-(g-1)u} {m_i \choose u}
} < 2 c^2 |\FF_{i, Y}|^g. 
\eeq 
Delete $Y$ not meeting this from $\YY_i$, and we have  
$|\YY_i| > \lp 1 - c^{-1} \rp {n_* \choose l}$, and \refeq{eq2Step1} $\wedge$ \refeq{eqStep2_2} for each remaining $Y \in \YY_i$.

\medskip

{\bf Step 3.~}{\em Show a similar fact with $\Pi(j+1)$-iii).} 
In our arguments of Step 2, replace the occurrences of the interval $[r_*, m_i ]$ by $[m_i]$, and those of the condition $\Pi(j)$-ii) by 
\[
\sum_{u \in [m_i] \atop T \in {Z_{i, j} \choose u}}~ \frac{|\FF_i[T]|^g}{\lp \frac{b_j}{12} \rp^{-(g-1)u} {m_i \choose u}} 
< |\FF_i|^g, 
\]
which is immediate from $\Pi(j)$-iii). This updates $\TT$ into 
$
\bigcup_{u=[m_i]} {Z \choose u}
$, and $b$ into $b_j \big/ 12$. 
Follow the step indentically except for the changes. After its completion, we have 
\beeq{eqStep3}
\sum_{
u \in [m_i] 
\atop 
T \in {Z \choose u}
}~ 
\frac{|\FF_{i, Y}[T]|^g  
}
{
\lp \frac{b_j}{12} \rp^{-(g-1)u} {m_i \choose u}
} < 2 c^2 |\FF_{i, Y}|^g, 
\eeq
for  more than $\lp 1-  c^{-1} \rp {n_* \choose l}$ sets $Y \in {X_* \choose l}$. Eliminate the other $Y$ from $\YY_i$ to have $|\YY_i| > \lp 1 - 3 c^{-1} \rp {n_* \choose l}$, and \refeq{eq2Step1} $\wedge$ \refeq{eqStep2_2} $\wedge$ \refeq{eqStep3} for $Y \in \YY_i$.

\medskip

{\bf Step 4.~}{\em Find a tuple $(Y_1, Y_2, \ldots, Y_k)$ of mutual disjoint $k$ sets $Y_i \in \YY_i$.} Extend $\YY_i$ to the family $\YY'_i$ of $l'$-sets where $l'= l \cdot 2^{\lc  \lg  \lg k \rc}$. This means 
\[
\YY'_i = Ext \lp \YY_i, l' \rp, \eqand  l' < \frac{n_* \sqrt \ep}{k}.  
\]

By Remark D) of Sec.\ 2.2, 
\[
\kappa \lbr {X_* \choose l'}  -  \YY'_i  \rbr
\ge   \kappa \lbr {X_* \choose l}  -  \YY_i  \rbr \frac{l'}{l} >  2^{\lg \lg k}  \ln \frac{c}{3}. 
\] 
Its RHS is no less than
\[
\frac{\ln c/3}{\ln 2}  \ln k  > \ln \frac{k}{\delta}, \eqwhere 
\delta:= k^{-\frac{\ln c/3}{\ln 2} +1}. 
\]
As $k \ge 2$ WLOG, $\delta$ dimishes indefinitely as $c > 2^{1/\ep}$ grows larger. With the $\delta$ we have 
\beeq{eqStep4}
|\YY'_i| > \lp 1 - \frac{\delta}{k} \rp {n_* \choose l'}. 
\eeq

Now consider all $i \in [k]$. There exists a tuple $(Y'_1, Y'_2, \ldots, Y'_k)$ of mutually disjoint $k$ sets $Y'_i \in \YY'_i$. For, there are $\prod_{i=1}^k {n_* -(i-1) l' \choose  l'}$ tuples $(Y'_1, Y'_2, \ldots, Y'_k)$ of mutually disjoint sets $Y'_i \in {X_* \choose l'}$.  By \refeq{eqStep4}, less than $\delta k^{-1}$ of them fail to satisfy $Y'_i \in \YY'_i$ for each particular $i$. Hence, the remaining more than $1- \delta$ of them each meet $Y'_i \in \YY'_i$ for all $i$.

Fixing such a $(Y'_1, Y'_2, \ldots, Y'_k)$, we obtain a set tuple $(Y_1, Y_2, \ldots, Y_k)$ of mutually disjoint $l$-sets $Y_i \in \YY_i \cap {Y'_i \choose l}$.  Update each $\FF_i$ by $\FF_i \leftarrow \FF_{i, Y_i}$. Then the following conditions hold for each $i$ due to \refeq{eq2Step1} $\wedge$ \refeq{eqStep2_2} $\wedge$ \refeq{eqStep3} with $Y \leftarrow Y_i$:

\beeqn
\textrm{i)}
&& \label{eqStep5} 
|\hspace*{0.5mm} \textrm{current}~\FF_i \hspace*{0.5mm}| ~>~ 
\lp 1 - \frac{1}{cg} \rp
\frac{{l  \choose m_*}}{{n_* \choose m_*}}
~|\hspace*{0.5mm}\FF_i~\textrm{before Step 1}\hspace*{0.5mm}|, 
\\ \textrm{ii)} && \nonumber 
\sum_{u \in [r_*, m_i] \atop T \in {Z \choose u}} 
\hspace*{3mm}
\frac{|\FF_i[T]|^g}
{
b_{j+1}^{-(g-1)u} {m_i \choose u} 
}
< \ep^2 |\FF_i|^g, 
\\ \textrm{iii)} && \nonumber 
\sum_{u \in [m_i] \atop T \in {Z \choose u}} 
\hspace*{3mm}
\frac{|\FF_i[T]|^g}
{
b_{*, j}^{-(g-1)u} {m_i \choose u} 
}
< |\FF_i|^g, 
\eqwhere b_{*, j} := \frac{b_j}{24 c^2}. 
\\ \textrm{iv)} && \nonumber 
\textrm{$\FF_i$ is element-wise disjoint with the other $\FF_{i'}$ within $\bigcup_{j'=1}^{j+1} X_{j'}$,}
\eeqn 
Note that $b$ has been replaced by $b_{j+1} = \lp 1- r^{-1} \rp b$, which changes $2 c^2 |\FF_{i, Y_i}|^g$ of \refeq{eqStep2_2} into $\ep^2 |\FF_i|^g$ of ii); it is due to $\lp 1- r^{-1} \rp^u < c^{-1}$ for $u \in [r_*, m_i]$.

\medskip

{\bf Step 5.~}{\em Finalize the induction step.} Execute the algorithm {\sc GrowCores} in \reffig{fig1} to further update each $C_i$ and $\FF_i$. Observe extra remarks.

\begin{figure} 
\begin{small}
\be{}
\item {\sc GrowCores}
\item ~{\em Inputs:} 
\\ \hspace*{1mm}
$C_i$ given by $\Pi(j)$ and $\FF_i$ by Step 4 for $i \in [k]$. 
\item ~{\em Outputs:}  
\\ \hspace*{1mm}
Updated $C_i$ and $\FF_i$ that satisfy $\Pi(j+1)$-i) to iv). 
\ee 
\be{\hspace*{2.5mm} 1.} 
\item \bfor $i=1$ \bto $k$ \bdo
\be{{1-}1.} 
\item perform substitutions into the tentative variables as follows: 
\(
&& 
X \leftarrow X - C_i; \quad 
\FF \leftarrow \lb U - C_i~:~ U \in \FF_i \rb; \quad 
m \leftarrow m_i; \quad 
\\ && 
u_* \leftarrow r_*; \quad 
b \leftarrow b_{j+1}; \quad 
h \leftarrow \ep; \quad 
Z \leftarrow Z; 
\)
\item apply the push-up lemma (\reflm{pushup} in Sec.\ 2.3) to the tentative variables with \refeq{eqStep5}-ii) to find a set $S$; 
\item $S_i \leftarrow S$; \hspace*{1mm}
$C_i \leftarrow C_i \cup S_i$; \hspace*{1mm}
$\FF_i \leftarrow \FF_i[C_i]$; \hspace*{1mm}
\item \bfor $i' \in [k] - \lb i \rb$ \bdoo $\FF_{i'} \leftarrow \FF_{i'} - \FF_{i'}(C_i)$; 
\ee 
\ee 
\end{small}
\caption{Algorithm {\sc GrowCores} to Further Update $C_i$ and $\FF_i$}
\label{fig1} 
\end{figure}

\be{A)} \setcounter{enumi}{5} 
\item Before the algorithm starts, every $\FF_i$ satisfies the $\Gamma \lp c^6 g^6 k_\dag  \sqrt m \rp$-condition in $Z_{i, j}$, confirmed similarly to Step 1 with \refeq{eqStep5}-ii). 

\item After Step 1-3 of {\sc GrowCores} with $i=1$, the family $\FF_1$ satisfies 
\( 
\textrm{i)}
&& 
\sum_{u \in [r_* -|S_i|, m_i- |S_i|] \atop T \in {Z-S_i \choose u} 
}
\hspace*{2mm}
\frac{|\FF_i[S_i \cup T]|^g}{b_{j+1}^{-(g-1)u} {m_i-|S_i| \choose u}}
< \ep |\FF_i[S_i]|^g. 
\\ \textrm{ii)} && 
\sum_{u \in [r_* - |S_i|] \atop T \in {Z-S_i \choose u} 
}
\hspace*{2mm}
\frac{|\FF_i[S_i \cup T] |^g}{\lp \frac{b_{j+1}}{6} \rp^{-(g-1)u} {m_i-|S_i| \choose u}}
< |\FF_i[S_i] |^g. 
\\ \textrm{iii)} && 
|\FF_i[S_i] | >  \ep^{\frac{1}{g-1}} b_{j+1}^{-|S_i|} |\FF_i |, 
\)
due to the push-up lemma referred to by Step 1-2. By ii), $\FF_i[S_i]$ satisfies the $\Gamma \lp c^6 g^6 k_\dag \sqrt m \rp$-condition in $Z_{i, j}-S_i$. 

\item 
By Step 1-4 with $i=1$, each $\FF_{i'}$ is reduced by less than its $(c^2 g k)^{-1}$, since F) with $|S_1| < r_* < c^6 g^6 k_\dag \sqrt m \big/ (c^2 g k)$. 
\item $\FF_2$ satisfies the $\Gamma\lbr c^6 g^6 k_\dag  \sqrt m,~ 1 - (c g k)^{-1} \rbr$-condition in $Z_{i, j}$ before Step 1-3 with $i=2$, and the $\Gamma\lp c^6 g^6 k_\dag \sqrt m \rp$-condition in $Z_{i, j} - S_i$ after it. 
\item Likewise, the other 
$\FF_i$ satisfy the $\Gamma\lbr c^6 g^6 k_\dag \sqrt m,~ 1 - (i-1) (c g k)^{-1} \rbr$-condition in $Z_{i, j}$ before Step 1-3 with the loop index $i$, and the $\Gamma\lp c^6 g^6 k_\dag \sqrt m \rp$-condition in $Z_{i, j} - S_i$ after it. When the algorithm terminates, every $\FF_i$ satisfies the \\ $\Gamma\lbr c^6 g^6 k_\dag \sqrt m,~ 1 - (c g) ^{-1} \rbr$-condition because $\FF_i$ possibly reduces by less than its $(cg)^{-1}$ since Step 1-3 with the index $i$. 
\item 
Consequently, the final $C_i$ and $\FF_i$ from {\sc GrowCores} satisfy:  
\( 
\textrm{i)}
&& 
\sum_{u \in [r_* -|S_i|, m_i- |S_i|] \atop T \in {Z-S_i \choose u} 
}
\hspace*{2mm}
\frac{|\FF_i[T]|^g}{b_{j+1}^{-(g-1)u} {m_i-|S_i| \choose u}}
< 2 \ep |\FF_i|^g. 
\\ \textrm{ii)} && 
\sum_{u \in [r_* - |S_i|] \atop T \in {Z-S_i \choose u} 
}
\hspace*{2mm}
\frac{|\FF_i[T] |^g}{\lp \frac{b_{j+1}}{6} \rp^{-(g-1)u} {m_i-|S_i| \choose u}}
< (1+\ep) |\FF_i |^g. 
\qed 
\)
\ee 

\medskip 

We select the obtained $C_i$ and $\FF_i$ as the ones for $\Pi(j+1)$, for which we verify the four conditions: $\Pi(j+1)$-i) holds since $\Pi(j)$-i) and $|S_i| < r_*$ from the push-up lemma. The conditions $\Pi(j+1)$-ii) and iii) hold due to K) while iv) is from \refeq{eqStep5}-iv).

We have proven the induction step $\Pi(j) \Rightarrow \Pi(j+1)$ when $m_*>0$ for all $i$. If $m_*=0$ for an $i$ before Step 1, choose the given $C_i$ and $\FF_i$ as the final ones. It is straightforward to check that they satisfy the four conditions of $\Pi(j+1)$ as well. This completes the induction step. 

\medskip

We now have $\Pi(r-1)$ implying the following four: 
\bdash 
\item $|C_i| < 2m/r$ by its i). 
\item The $\Gamma( \hat b)$-condition of each $\FF_i$ in $Z_{i, r-1}$ where 
$\hat b = c^6 g^6 k_\dag m_*$ with  
\\ $m_*= q_{i, r}- |C_i \cap X_r|$ that is positive WLOG. It is true by the same arguments as Step 1 above. 
\item $|\FF_i| > \hat b^{m_*}$ implied by the $\Gamma$-condition. 
\item The element-wise disjointness of $\FF_i$ within $\bigcup_{j'=1}^{r-1} X_{j'}$. 
\edash 

For each $i$, perform the same as Step 1 to find a family $\YY_i \subset {X_* \choose l}$ such that
$
|\YY_i| > {n_* \choose l} \lp 1  - \frac{1}{cg} \rp
$ and \refeq{eq2Step1}. It means 
\(
&& 
\frac{2|\FF_{i, Y}|}{|\FF_i|} > \frac{{l \choose m_*}}{{n_* \choose m_*}} 
= \prod_{u=0}^{m_*-1} \frac{l-u}{n_*-u} > \lp \frac{\ep}{2 k_\dag} \rp^{m_*}, 
\sothat 
|\FF_{i, Y}| > (c^5 g^6 m_*)^{m_*}, 
\)
for $Y \in \YY_i$.

Then the same construction as Step 4 confirms the existence of a set tuple $(Y_1, Y_2, \ldots, Y_k)$ such that $Y_i \in \YY_i$ that are mutually disjoint.  By this we have obtained $k$ nonempty subfamilies $\FF_{i, Y_i}$ that are element-wise disjoint within $X$. As they include $k$ mutually disjoint sets, we have completed our proof of \refth{KDC}. 


\appendix  
\section{Proof of Some Statements in Sections 1 and 2}
\setcounter{section}{1}

\subsection{Relations in $\lp 2^X \rp^g$ and $\Gamma$-Conditions}  
We prove Lemmas \ref{gthMark}, \ref{lmConversion} and \ref{pushup} in this subsection.

\medskip 

{\noindent \bf \reflm{gthMark}.}  {\em ($g^{th}$ Mark Lemma) \cite{3SF} Let 
\be{i)} 
\item $X$ be the universal set weighted by $w: \lp 2^X \rp^g \rightarrow \R_{\ge 0}$ for some $g \in \Z_{\ge 2}$, 
\item $l \in [n]$, $m \in [l-1]$, and $h, \gamma \in \R_{>0}$, such that 
$\gamma$ and $\frac{l}{gm}$ are both sufficiently large, 
\item and $\FF \subset {X \choose m}$ satisfy the $\Gamma_g  \lp \frac{4 \gamma n}{l},~h \rp$-condition on the norm $\| \cdot \|$ induced by $w$.  
\ee 
Then 
\[ 
\| \DD_g  \| < 
\frac{\lp 1 + \frac{h}{\gamma} \rp {n \choose l}{l \choose m}^g}{{n \choose m}^g} 
\| \FF \|^g.
\] 
}
\begin{proof}
Given the objects, write 
\[
b = \frac{4 \gamma n}{l}, \eqand 
v(x) = {x \choose m}^{-g+1}  \hspace*{3mm} \prod_{p=1}^{g-1} {x-pm \choose m},  
\]
for $x \in [n] - [g m]$.  Observe the five remarks. 
\be{a)} 
\item $\PP_j = \emptyset$ if $j > (g-1)m$, so 
\(
\|  \DD  \|  &=& 
\sum_{Y \in{X \choose l}} ~\sum_{\bs{U} \in {Y \choose m}^g} 
w \lp \bs{U} \rp
= 
\sum_{Y \in {X \choose l}} \left\| {Y \choose m} \right\|^g 
\nexteqline
\sum_{j=0}^{(g-1)m} \| \PP_j \| {n- gm+j \choose l - gm+j}. 
\) 
Due to the $\Gamma_g \lp b, h\rp$-condition of $\FF$ and $\| \FF \|^g = {n \choose m}^g e^{-g \kappa \lp \FF \rp}$, 
\[
\|  \DD  \| - \| \PP_0 \| {n- gm \choose l - gm}< 
h  \sum_{j=1}^{(g-1)m} b^{-j}  {n \choose m}^g e^{- g\kappa\lp \FF \rp}
{n-gm +j \choose l-gm+j}. 
\]

\item  
$
{v(l) \le v(n)} 
$,  
since 
\(
&&
\prod_{i=0}^{m-1} \frac{1-\frac{pm}{l-i}}{1-\frac{pm}{n-i}} \le 1, \quad 
\textrm{for~} p \in [g-1], 
\sothat 
\frac{{l- pm \choose m}{l \choose m}^{-1}}{{n-pm \choose m}{n \choose m}^{-1}} \le 1, 
\quad \Rightarrow \quad 
\frac{v(l)}{v(n)} \le 1. 
\)

\item We have the identify 
\[
{n \choose m}^g {n-gm \choose l-gm}  
= \frac{v(l)}{v(n)} {n \choose l} {l \choose m}^g, 
\]
since its LHS equals 
\( 
&&
\frac{1}{v(n)}  {n-gm \choose l-gm} \prod_{p=0}^{g-1} {n-pm \choose m} 
\nexteqline 
\frac{{l-(g-1)m \choose m}}{v(n)}  {n-(g-1)m \choose l-(g-1)m} 
\prod_{p=0}^{g-2} {n-pm \choose m} 
\nexteqline 
\frac{{l-(g-1)m \choose m} {l-(g-2)m \choose m}}{v(n)}  {n-(g-2)m \choose l-(g-2)m} 
\prod_{p=0}^{g-3} {n-pm \choose m} 
\nexteqline 
\cdots =  
\frac{\prod_{p=0}^{g-1} {l-pm \choose m}}{v(n)} {n \choose l}
=
\frac{v(l)}{v(n)} {n \choose l} {l \choose m}^g, 
\) 
due to the identity ${x \choose z}{x-z \choose y-z}= {x \choose y}{y \choose z}$. 
By b), 
\[
{n \choose m}^g {n-gm \choose l-gm}  \le {n \choose l} {l \choose m}^g. 
\]
This inequality can be derived by counting the number of $g^{th}$ marks as well\footnote{Let $\FF= {X \choose m}$, and $w$ be the unit weight with respect to $\FF$, $i.e.$, $w \lp \bs{U} \rp = 1$ if all $U_i$ in $\bs{U}= (U_1, U_2, \ldots, U_g)$ are in $\FF$, and $w \lp \bs{U} \rp = 0$ otherwise. Then, 
\[
{n \choose l}{l \choose m}^g = \| \DD \| = \sum_{j=0}^{(g-1)m} \| \PP_j \| {n - gm + j \choose l - gm+j} \ge {n \choose m}^g {n - gm \choose l-gm}, 
\]
by a), $\sum_{j=0}^{(g-1)m} \| \PP_j \| = \| \FF \|^g = {n \choose m}^g$, and  ${n - gm + j \choose l - gm+j} \ge {n - gm \choose l-gm}$ for any $j$. 
}.

\item So, 
\( 
\| \PP_0 \| {n-gm \choose l-gm}
&\le&
\| \FF \|^g  {n-gm \choose l-gm} 
= {n \choose m}^g {n-gm \choose l-gm}   e^{- g \kappa \lp \FF \rp} 
\\ &\le & 
{n \choose l} {l \choose m}^g e^{- g \kappa \lp \FF \rp} . 
\)

\item For each $j \in [(g-1)m]$,  
\[
{n-gm +j \choose l-gm+j} 
= {n-gm \choose l-gm} \prod_{i=0}^{j-1} \frac{n-gm+j-i}{l-gm+j-i}  
<
\lp \frac{2n}{l} \rp^{j}
{n-gm \choose l-gm},  
\]
since $l$ and $n$ are both sufficiently larger than $gm$. 
\ee

\medskip

By the remarks, it suffices to show that 
\beeq{eqTildeG} 
\| \DD \| < \lbr  1  
+ h \sum_{j=1}^{(g-1)m} \lp 2 \gamma \rp^{-j}  
\rbr 
{n \choose l}{l \choose m}^g e^{-g \kappa \lp \FF \rp}, 
\eeq 
as its RHS is less than 
$
\frac{\lp 1 + \frac{h}{\gamma} \rp {n \choose l}{l \choose m}^g}{{n \choose m}^g}
\sum_{\bs{U} \in \FF^g} w \lp \bs{U}  \rp$. We see from a), c) and e) that 
\(
&&
\sum_{j=1}^{(g-1)m} \| \PP_j  \|  {n-gm +j \choose l-gm+j} 
\\ &<& \nonumber 
h  \sum_{j=1}^{(g-1)m} b^{-j}   {n \choose m}^g e^{- g\kappa\lp \FF \rp}
{n-gm +j \choose l-gm+j} 
\\ &<& \nonumber 
h e^{-g \kappa \lp \FF \rp} 
{n \choose m}^g  {n-gm \choose l-gm} 
\sum_{j=1}^{(g-1)m} \lp 2 \gamma  \rp^{-j} 
\\ &\le& \nonumber 
h e^{-g \kappa \lp \FF \rp} 
{n \choose l}{l \choose m}^g 
\hspace*{2mm}
\sum_{j=1}^{(g-1)m} \lp 2 \gamma \rp^{-j}.  
\) 
Also by d), this confirms \refeq{eqTildeG} completing the proof of the lemma. 
\end{proof}

\medskip 

{\noindent\bf \reflm{lmConversion}.} 
{\em (Conversion Lemma) 
Let $X$ be weighted by $w: \lp 2^X \rp^g \rightarrow \R_{\ge 0}$ inducing $\| \cdot \|$, primitively with some norm $\| \cdot \|_*$. If a family $\FF \subset {X \choose m}$ satisfies the $\Gamma(b, h)$-condition on $\| \cdot \|_*$ for some $b, h \in \R_{\ge 1}$, it satisfies the $\Gamma_g \lbr \frac{b}{2^{g-2} (g-1) m},~ h^{g-1} \rbr$-condition on $\| \cdot \|$. 
}
\begin{proof}
Write $b_{g'} = 2^{-g' + 2} b$ for $g' \in [2, g]$ to start our proof. 
It suffices to confirm 
\beeq{eqSec14}
\| \PP_{j, g} \| < h^{g-1} b_g^{-j} {(g-1)m \choose j} \| \FF \|^g, \quad \forall j \in [(g-1)m]. 
\eeq 
Define 
\[ 
w_{g'}: \lp 2^{X}\rp^{g'} \rightarrow \R_{\ge 0}, \quad 
(U_1, U_2, \ldots, U_{g'}) \mapsto 
\prod_{i=1}^{g'} w_*(U_i), 
\] 
for each $g'$ inducing the norm $\| \cdot \|_{g'}$, where $w_*: 2^X \rightarrow \R_{\ge 0}$ is the primitive weight that induces the given $\| \cdot \|_*$. We verify 
\beeq{eqSec14_2}
\| \PP_{j, g'} \|_{g'} < h^{g'-1} b_{g'}^{-j} {(g'-1)m \choose j} 
\| \FF \|_{g'}^{g'}, 
\eeq 
for every $g'$ and $j \in [(g'-1)m]$. The case $g' = g$ proves \refeq{eqSec14}. 

\medskip 

Proof of \refeq{eqSec14_2} by induction on $g'$. Denote $g'$ by $g$ for simplicity. Fix each $j \in [(g-1) m]$ for the basis $g=2$ to have 
\(
\| \PP_{j, 2} \|_2 &=& \sum_{\bs{U} \in \PP_{j, 2}} w_2 \lp  \bs{U} \rp  
= 
\sum_{(U_1, U_2) \in \PP_{j, 2}} \|\FF[U_1] \|_*~\|\FF[U_2]\|_*
\\ &=& 
\sum_{U_1 \in \FF} \|\FF[U_1] \|_*
\sum_{U_2 \in \FF \atop \textrm{with~} |U_1 \cap U_2|=j} \|\FF[U_2]\|_*
\\ &<& 
h  b_2^{-j} 
{m \choose j} \|\FF \|_2^2
, 
\)
because: 
\bdash
\item By Remark A) of Sec.\ 2.2, $\sum_{U_1 \in \FF} \|\FF[U_1] \|_* = \|\FF \|_* = \|\FF \|_2$ since $\| \cdot \|_2$ is primitive with $\| \cdot \|_*$. 
\item For each $U_1 \in \FF$ and $S \in {U_1 \choose j}$, the sum of $\|\FF[U_2]\|_*$ over $U_2 \in \FF[S]$ is $\| \FF[S] \|_* < h b^{-j} \| \FF \|_*$ by the given $\Gamma(b. h)$-condition of $\FF$. 
\edash 
This proves the basis.

Assume true for $g-1$ and prove true for $g$. For each $j \in [(g-1)m]$, 
\(
&& 
\sum_{k \in [j] \atop \bs{U} \in \PP_{k, g-1}}  w_{g-1} (\bs{U})  
\sum_{S \in {Union(\bs{U}) \choose j-k}} w_*(S) 
\\ &<& 
\sum_{k=1}^{j}
h^{g-2} b_{g-1}^{-k} {(g-2)m \choose k} \| \FF \|_{g-1}^{g-1} 
\cdot 
{(g-1) m -k  \choose j-k} h b^{-j+k} \| \FF \|_*,  
\)
since: 
\bdash
\item $\PP_{k, g-1}=\emptyset$ if $k > (g-2) m$. 
\item By induction hypothesis, 
$ 
\sum_{\bs{U} \in \PP_{k, g-1} 
} 
w_{g-1} \lp \bs{U} \rp
< 
h^{g-2} b_{g-1}^{-k} {(g-2)m \choose k} \| \FF \|_{g-1}^{g-1} 
$ 
for $k \in [(g-2)m]$. 
\item $|Union(\bs{U})| = (g-1)m-k$ for each $\bs{U} \in \PP_{k, g-1}$ by definition. 
\edash

Also noting $
{(g-2)m \choose k} {(g-1) m -k  \choose j- k} 
< 
{(g-1)m \choose k} {(g-1) m -k  \choose j- k} 
= 
{(g-1)m \choose j}{j \choose k}
$ for every $k$, we get 
\(
\sum_{k \in [j] \atop \bs{U} \in \PP_{k, g-1}}  w_{g-1} (\bs{U})  
\sum_{S \in {Union(\bs{U}) \choose j-k}} w_*(S) 
&<& 
h^{g-1} b_{g-1}^{-j} 
{(g-1) m \choose j} 
\| \FF \|_g^g 
\sum_{k=1}^{j}
{j \choose k}
\\ &=& 
h^{g-1} \lp 2 b_g \rp^{-j} (2^j - 1) 
{(g-1) m \choose j} 
\| \FF \|_g^g. 
\)

For $k=0$, 
\(
\sum_{\bs{U} \in \PP_{0, g-1}}  w_{g-1} (\bs{U})  \sum_{S \in {Union(\bs{U}) \choose j}} w_*(S) 
&<& \| \FF \|_*^{g - 1} {(g-1)m \choose j} h b^{-j} \| \FF \|_*
\\ &\le& h^{g-1} \lp 2 b_g \rp^{-j} 
{(g-1) m \choose j}  \| \FF \|_g^g, 
\)
by $\sum_{\bs{U} \in \PP_{0, g-1}}  w_{g-1} (\bs{U}) \le  \| \FF \|_{g-1}^{g - 1}$, $h \ge 1$ and $g-1 \ge 2$. Therefore, 
\(
\| \PP_{j, g} \|_g  
&\le & 
\sum_{k \in [0, j] \atop \bs{U} \in \PP_{k, g-1}}  w_{g-1} (\bs{U})  
\sum_{S \in {Union(\bs{U}) \choose j-k}} w_*(S) 
\\ &<& 
h^{g-1} b_g^{-j} 
{(g-1) m \choose j} 
\| \FF \|_g^g, 
\)
completing the induction step.

This confirms \refeq{eqSec14_2} and the desired $\Gamma_g$-condition of $\FF$, proving the conversion lemma. 
\end{proof}

\medskip 

{\noindent\bf \reflm{pushup}.}
{\em (Push-Up Lemma) 
In the universal set $X$ weighted with a primitive norm $\| \cdot \|$, let $\FF \subset {X \choose m}$ satisfy 
\[
\sum_{u \in [u_*, m] \atop S \in {Z \choose u}
}
\frac{\|\FF[S]\|^g}{b^{-(g-1)u} {m \choose u}}
< h^2 \|\FF \|^g, 
\]
for some set $Z \subset X$, integers $g \in \Z_{\ge 2}$ and $u_* \in [0, m]$, and real numbers $b \in \R_{>1}$ and $h \in (0, 1/4]$. Then there exists $S \in {Z \choose v}$ for some $v \in [0, u_*)$ satisfying the three conditions. 
\( 
1) && 
\sum_{u \in [u_* -v, m-v] \atop 
T \in {Z-S \choose u} 
}
\frac{\|\FF[S \cup T]\|^g}{b^{-(g-1)u} {m-v \choose u}}
< h \|\FF[S] \|^g. 
\\ 2) && 
\sum_{u \in [u_* - v] \atop 
T \in {Z-S \choose u} 
}
\frac{\|\FF[S \cup T] \|^g}{\lp \frac{b}{6} \rp^{-(g-1)u} {m-v \choose u}}
< \|\FF[S] \|^g. 
\\ 3) && 
\|\FF[S] \| >  h^{\frac{1}{g-1}} b^{-v} \|\FF\|. 
\)
} 
\begin{proof}
First find the maximum $v$ such that 
\beeq{eqPushup}
\sum_{ 
S \in {Z \choose v} 
}
\frac{\|\FF[S]\|^g}{b^{-(g-1)v} {m \choose v}}
\ge 4 h \|\FF \|^g. 
\eeq 
There does exist such $v \in [0, u_*)$ because the inequality holds when $v=0$, and would not if $v=u_*$.

\medskip 

We show a property here. 

\medskip 

{\em Claim 1.}
\[
\sum_{S \in {Z \choose v} \atop \textrm{with $\neg$ 3)}}
\frac{\|\FF[S]\|^g}{b^{-(g-1)v} {m \choose v}} \le h \| \FF\|^g. 
\]
\\{\em Proof.}  
Consider 
\[
\sum_{S \in \SS} \|\FF[S]\|^g, \eqwhere \SS = \lb S~:~ S\in {Z \choose v} ~\textrm{with $\neg$ 3)}~\rb, 
\] 
noting that the summand $\|\FF[S]\|^g$ equals the sum of $\prod_{i=1}^g \|  U_i  \|$ for all $(U_1, U_2, \ldots, U_g) \in \|\FF[S]\|^g$, since $\| \cdot \|$ is primitive.  

Fix each $U \in \FF$ to observe that 
\[ 
\| U \| \sum_{ 
S \in {U \choose v} \cap \SS} \| \FF[S] \|^{g-1}
\le 
\| U \|~{m \choose v} h b^{-(g-1)v} \| \FF\|^{g-1}, 
\]  
because$\| \FF[S] \|^{g-1}$ is no more than $h b^{-(g-1)v} \| \FF\|^{g-1}$ due to the condition $\neg$ 3) of $S$. 
This means 
\[
\sum_{S \in \SS} \|\FF[S]\|^g \le  h b^{-(g-1)v} {m \choose v} \| \FF\|^g, 
\]
which the claim follows.

\medskip 

By Claim 1 and \refeq{eqPushup}, 
$
\sum_{S \in {Z \choose v} \atop \textrm{with $\neg$ 1) or $\neg$ 2)}}
\frac{\|\FF[S]\|^g}{b^{-(g-1)v} {m \choose v}} \ge 3h \| \FF\|^g
$, 
if there were no $S \in {Z \choose v}$ such that 1) $\wedge$ 2) $\wedge$ 3). It means one of the following two would be true: 
\(
- && 
\sum_{S \in {Z \choose v} \atop \textrm{with~} \neg 1)}  
\frac{\|\FF[S] \|^g}{b^{-(g-1)v} {m \choose v}}
\ge  h \| \FF \|^g; 
\\ - && 
\sum_{S \in {Z \choose v} \atop \textrm{with~} \neg 2)\textrm{-}u}  
\frac{\| \FF[S] \|^g}{b^{-(g-1)v} {m \choose v}}
\ge  \frac{4 h}{3^u} |\FF|^g, 
\quad 
\textrm{for some~} u \in [u_*-v], 
\)
where 2)-$u$ means 
\[
\sum_{
T \in {Z-S \choose u} 
}
\frac{\| \FF[S \cup T] \|^g}{\lp \frac{b}{6} \rp^{-(g-1)u} {m-v \choose u}}
< 2^{-u} \|\FF[S] \|^g. 
\] 
Call the two cases {\em Cases 1} and {\em 2}, respectively. 

\medskip 

We show a contradiction in Case 1 from 
\beeqn  
&& \label{eq3Pushup} 
\sum_{S \in {Z \choose v} \atop \textrm{with~} \neg 1)} 
\sum_{u \in [u_* - v, m - v] \atop 
T \in {Z-S \choose u} 
}
\frac{|\FF[S \cup T]|^g}{b^{-(g-1)(v+u)} {m \choose v+u} {v + u \choose v}} 
\nexteqline \nonumber 
\sum_{S \in {Z \choose v} \atop \textrm{with~} \neg 1)} 
\frac{1}{b^{-(g-1)v} {m \choose v}}
\sum_{ u \in [u_*-v, m-v] \atop 
T \in {Z-S \choose u} 
}
\frac{\| \FF[S \cup T] \|^g}{b^{-(g-1)u} {m-v \choose u}} 
\\ &\ge& \nonumber 
h \sum_{S \in {Z \choose v} \atop \textrm{with~} \neg 1)} 
\frac{\| \FF[S] \|^g}{b^{-(g-1)v} {m \choose v}} 
\\ &\ge& \nonumber 
h^2  |\FF|^g.  
\eeqn 
The second line is true due to ${m \choose v}{m- v \choose (v+u)-v}= {m \choose v+u} {v+u \choose v}$, and the third line to each $S$ meeting $\neg$ 1). 
The relation means: 

\medskip 

{\em Claim 2.} 
\beeq{eq4Pushup} 
\sum_{
v' \in [u_*, m] \atop 
S' \in {Z \choose v'}} 
\frac{\| \FF[S'] \|^g}{b^{-(g-1)v'} {m \choose v'}}  
\ge h^2 \| \FF \|^g.  
\eeq 
\\{\em Proof.} 
Let 
\[
\gamma_{v'} =  \| \FF \|^{-g}  \sum_{S' \in {Z \choose v'}} \frac{\| \FF[S'] \|^g}{b^{-(g-1) v'} {m \choose v'}}, 
\]
for each $v' \in [u_*, m]$. 
It satisfies 
\[
\sum_{S \in {Z \choose v},~T \in {Z-S \choose v'- v}} 
\| \FF[S \cup T] \|^g 
= 
\sum_{S' \in {Z \choose v'}} \| \FF[S'] \|^g  {v' \choose v}
= 
\gamma_{v'} \| \FF \|^g b^{-(g-1) v'} {m \choose v'} {v' \choose v}
,  
\]
since there are ${v' \choose v}$ pairs $(S, T)$ such that $S \cup T$ equals each given $S' \in {Z \choose v'}$. 
Summing up 
\[
\sum_{S \in {Z \choose v},~T \in {Z-S \choose v'- v}} 
\frac{\| \FF[S \cup T] \|^g }{b^{-(g-1)v'} {m \choose v'}{v' \choose v}}  
= \gamma_{v'} \| \FF \|^g 
. 
\]
for all $v'$, we see $\neg$ \refeq{eq4Pushup} $\Rightarrow$ 
$\neg$ \refeq{eq3Pushup}. So \refeq{eq4Pushup}, proving Claim 2.

\medskip 

\noindent 
This contradicts the given collective $\Gamma$-condition, so Case 1 is impossible to occur.

\medskip 

Given $u \in [u_*-1]$ in Case 2, similarly find 
\(
&& 
\sum_{S \in {Z \choose v} \atop \textrm{with~} \neg 2)\textrm{-}u} 
\sum_{ 
T \in {Z-S \choose u} 
}
\frac{\| \FF[S \cup T] \|^g}{b^{-(g-1)(v+u)} {m \choose v+u} {v + u \choose v}} 
\\ &\ge&  
3^{(g-1)u}
\sum_{S \in {Z \choose v} \atop \textrm{with~} \neg 2)\textrm{-}u} 
\frac{\| \FF[S] \|^g}{b^{-(g-1) v} {m \choose v}} 
> 4 h \lp \frac{3^{g-1}}{3} \rp^u \| \FF \|^g, 
\)
meaning 
\[
\sum_{ 
S' \in {Z \choose v + u}} 
\frac{\| \FF[S'] \|^g}{b^{-(g-1) (v+u)} {m \choose v+u}}  
> 4 h \| \FF \|^g.  
\]
It is against the maximality of $v$ such that \refeq{eqPushup}. 

By this contradiction, we have shown that an $S \in {Z \choose v}$ satisfies 1) $\wedge$ 2) $\wedge$ 3), confirming the push-up lemma. 
\end{proof}


\subsection{Proof of \refth{EGT42}}  

This subsection proves: 

\medskip 

{\noindent\bf \refth{EGT42}.} \cite{3SF}
{\em Let $X$ be primitively weighted to induce the norm $\| \cdot \|$. 
For every sufficiently small $\ep \in (0, 1)$, and $\FF \subset {X \choose m}$ satisfying the $\Gamma_2 \lp \frac{4 \gamma n}{l},~1  \rp$-condition on $\| \cdot \|$ for some $l \in [n]$, $m \in [l]$, and 
$
\gamma \in \lbr \ep^{-2},~ l m^{-1} \rbr 
$, there are $\lc {n \choose l} \lp 1-  \ep \rp \rc$ sets $Y \in {X \choose l}$ such that
\[
\lp 1 - \sqrt{\frac{2}{\ep \gamma}} \rp
\frac{{l \choose m}}{{n \choose m}}
\left\| \FF \right\|
<
\left\| {Y \choose m} \right\|
<
\lp 1 + \sqrt{\frac{2}{ \ep \gamma}}  \rp
\frac{{l \choose m}}{{n \choose m}}
\left\| \FF \right\|. \qed 
\]
} 

\medskip 

Start our proof given such an $\FF$, $l$, and $\gamma$. We will refer to \refeq{eqTildeG} shown in the proof of the $g^{th}$ mark lemma, which holds true for the $\FF$, $l$, $\gamma$, $g=2$ and $h=1$ here. 
We also show another property on $\| \DD \|$.

\begin{lemma} \label{WeightBounds}
Let $l \in [n]$, $m \in [l]$, $t \in \R_{>0}$, and $\FF \subset {X \choose m}$ be weighted by $w$ so that
\(
0<
\| \DD \|
\le t {n \choose l} {l \choose m}^2 e^{- 2 \kappa \lp \FF \rp},
\)
and let $u, v \in \R_{>0}$ such that
\[
u<1, \hspace*{5mm}
u {n \choose l} \in \Z, \eqand
t < 1 + \frac{u (v-1)^2}{1-u}. 
\]
The two statements hold.
\begin{enumerate} [a)]
\item If $v \ge 1$, more than $(1-u) {n \choose l}$ sets $Y \in {X \choose l}$ satisfy
$
\left\|  {Y \choose m} \right\| <   v  {l \choose m} e^{-\kappa \lp \FF \rp}.
$
\item If $v \le 1$, more than $(1-u) {n \choose l}$ sets $Y \in {X \choose l}$ satisfy
$
\left\|  {Y \choose m} \right\| >   v  {l \choose m} e^{-\kappa \lp \FF \rp}.
$
\end{enumerate}
\end{lemma}
\begin{proof}
a): Put
\[
z = e^{-\kappa \lp \FF \rp}, \eqand
x_j = \left\| {Y_j \choose m} \right\|,
\]
where $Y_j$ is the $j$th $l$-set in ${X \choose l}$, and $z>0$ by $\| \DD \|>0$. Suppose to the contrary that
$x_j \ge v z {l \choose m}$ if $1 \le j \le u {n \choose l}$.

Noting
$
\sum_{1 \le j \le {n \choose l}} x_j = \| \MM \|  = z {n \choose l} {l \choose m}
$ from the above A),
let $y \in (0, 1)$ satisfy
\[
\sum_{1 \le j \le u {n \choose l}} x_j = y z {n \choose l}{l \choose m},
\]
so
$
y \ge uv.
$
Find
\(
&&
\sum_{1 \le j \le u {n \choose l}} x_j^2
\ge \lbr \frac{y z {n \choose l}{l \choose m}}{u {n \choose l}} \rbr^2  u {n \choose l}
= \frac{y^2 z^2}{u}  {n \choose l}{l \choose m}^2,
\\ \textrm{and} &&
\sum_{u {n \choose l} < j \le {n \choose l}} x_j^2
\ge \lbr \frac{(1-y) z {n \choose l}{l \choose m}}{(1-u) {n \choose l}} \rbr^2  (1-u) {n \choose l}
= \frac{(1-y)^2 z^2 }{1-u} {n \choose l}{l \choose m}^2,
\)
meaning
\beeqn
&& \label{eqWeightBounds}
\| \DD \| = \sum_{Y \in {n \choose l}} \left\| {Y \choose m} \right\|^2 \ge f z^2 {n \choose l}{l \choose m}^2,
\eqwhere
f = \frac{y^2}{u} + \frac{(1-y)^2}{1-u}.
\eeqn
From $y \ge u v  \ge u$,
\beeq{eq2WeightBounds}
f \ge uv^2 + \frac{(1-uv )^2}{1-u} = 1 + \frac{u (v-1)^2}{1-u}>t.
\eeq
This contradicts the given condition proving a).

\medskip

b): Suppose $x_j \le  v z {l \choose m}$ if $1 \le j \le u {n \choose l}$. Use the same $y$ and $f$ so
$y \le uv$ and \refeq{eqWeightBounds}. These also imply \refeq{eq2WeightBounds} producing the same contradiction. Thus b).
\end{proof}

Set 
\(
&& 
t= 
1  
+ \sum_{j=1}^{(g-1)m} \lp 2 \gamma \rp^{-j},
\quad 
u = \frac{\lf \frac{\ep}{2} {n \choose l} \rf}{{n \choose l}},
\\ && 
v = 1 + \frac{u}{\lp \frac{\ep}{2} \rp^{\frac{3}{2}}  \sqrt \gamma}, 
\eqand v' = 1- \frac{u}{\lp \frac{\ep}{2} \rp^{\frac{3}{2}}  \sqrt \gamma}. 
\)
Then $1 + \frac{u (v-1)^2}{1-u}= 1 + \frac{u (v'-1)^2}{1-u} 
>t$ since $\ep^{-2} \le \gamma \le l < {n \choose l}$ with the sufficiently small $\ep$.  Here $l<n$ is assumed as the theorem is trivially true if $l=n$. 
By \refeq{eqTildeG} as mentioned above, and \reflm{WeightBounds}, 
\[
v'   {l \choose m} e^{-\kappa \lp \FF \rp} <
\left\|  {Y \choose m} \right\|
<   v  {l \choose m} e^{-\kappa \lp \FF \rp},
\]
for some $\lp 1- 2 u \rp {n \choose l}$ sets $Y \in {X \choose l}$.
As $e^{-\kappa \lp \FF \rp} = \| \FF \|{n \choose m}^{-1}$, this means 
there are $\lc {n \choose l} \lp 1-  \ep \rp \rc$ sets $Y \in {X \choose l}$ such that
\[
\lp 1 - \sqrt{\frac{2}{\ep \gamma}}  \rp
\frac{{l \choose m}}{{n \choose m}}
\left\| \FF \right\|
<
\left\| {Y \choose m} \right\|
<
\lp 1 +\sqrt{\frac{2}{\ep \gamma}} \rp
\frac{{l \choose m}}{{n \choose m}}
\left\| \FF \right\|, 
\]
completing the proof of \refth{EGT42}.

\subsection{Confirmation of the Extension Generator Theorem}  
We now prove the extension generator theorem given in Sec.\ 1. 

\medskip 

{\noindent\bf \refth{EGT}.}
{\em (Extension Generator Theorem) \cite{3SF} 
There exists $\ep \in (0, 1)$ satisfying the following statement:
let $X$ be the universal set of cardinality $n$, $m \in [n-1]$, $l \in [n]-[m]$,  and
$
\lambda \in \lp 1, \frac{\ep l}{m^2} \rp.
$
For every nonempty family $\FF \subset {X \choose m}$, there exists an $\lp l, \lambda \rp$-extension generator $T$ of $\FF$ with  
$
|T| \le  \lbr \ln {n \choose m} - \ln |\FF| \rbr \Bigr/ \ln \frac{\ep l}{m^2 \lambda}.
$ \qed 
} 

\medskip

We note that \refco{corEGT42} implies:

\begin{corollary} \label{EGT4Cor}
Let $X$ be the universal set of cardinality $n$, $m \in [n-1]$, $l \in [n] - [m]$ and $\gamma \in \R_{>0}$ be sufficiently large not exceeding $\frac{l}{m^2}$. 
For any $\FF \subset {X \choose m}$ satisfying the $\Gamma \lp \frac{4  \gamma nm}{l} \rp$-condition, there are $\lc {n \choose l} \lp 1- \frac{2}{\sqrt[3] \gamma} \rp \rc$ sets $Y \in {X \choose l}$ such that
\[
\frac{{l \choose m} | \FF |}{{n \choose m}}
\lp 1 - \frac{1}{\sqrt[3] \gamma} \rp
<
\left| \FF \cap {Y \choose m} \right|
<
\frac{{l \choose m} | \FF |}{{n \choose m}}
\lp 1 + \frac{1}{\sqrt[3] \gamma} \rp. \qed
\]
\end{corollary}

\medskip 

We also use the following lemma. 

\medskip

\begin{lemma} \label{PhaseII} (Complement Sparsity Lemma) \cite{blog, 3SF} 
For $\FF \subset{X \choose m}$ such that $m \le \frac{n}{2}$,
\[
\kappa \lbr {X \choose 2m} - Ext \lp \FF, 2m \rp \rbr
\ge
2 \kappa \lbr {X \choose m} - \FF \rbr, 
\]
where the sparsity $\kappa$ is induced by the unit weight as in Remark C) of Sec.\ 2.2. 
\end{lemma}\begin{proof}
For each $S \in {X \choose m}-\FF$ and $j \in [0, m] \cap \Z$, let
\[
\FF_j = \lb T - S ~:~ T\in \FF,~ |T - S| = j \rb.
\]
There exists $j$ such that $\kappa \lp \FF_j \rp$ in the universal set $X- S$ is at most $\kappa \lp \FF \rp$ in X, otherwise
\[
|\FF| < \sum_{j\ge 0} {m \choose m-j}{n-m \choose j} e^{-\kappa \lp \FF \rp} = {n \choose m} e^{-\kappa \lp \FF \rp} = |\FF|.
\]
Taking $Ext \lp \FF_j, m \rp$ in $X-S$ referring to Remark E) of Sec.\ 2.2, we see
there are $\lc {n -m \choose m} e^{-\kappa \lp \FF \rp} \rc$ pairs $(S, U)$ such that $U \in {X-S \choose m}$ and $S \cup U \in Ext \lp \FF, 2m \rp$ for each $S \in {X \choose m}- \FF$.

Now consider all pairs $(S, U)$ such that $S$ and $U$ are disjoint $m$-sets, and $S \cup U \in Ext \lp \FF, 2m \rp$. Their total number is at least
${n \choose m}{n-m \choose m}={n \choose 2m}{2m \choose m}$
times
$
(1-z) + z (1-z) =1 - z^2
$
where $z=e^{-\kappa \lbr {X \choose m} - \FF \rbr}$.

As a $2m$-set produces at most ${2m \choose m}$ pairs $(S, U)$, there are at least $\lc {n \choose 2m}(1-z^2) \rc$ sets in $Ext \lp \FF, 2m \rp$.
The lemma follows.
\end{proof}

\medskip

We show \refth{EGT} with those two statements. Given sufficiently small $\ep$, $m, l, \lambda$ and $\FF$ as the statement, assume $|\FF| > b^m$ for some $b \in \R_{\ge 1}$. There exists $T \subset X$ such that $|T|<m$,  $|\FF[T]| \ge |\FF| b^{-|T|}$, and $|\FF[T \cup S]| < b^{m-|T \cup S|}|\FF[T]|$ for any nonempty $S \subset X - T$. We consider such $\FF[T]$ in the universal set $X-T$ in place of $\FF$ in our proof of \refth{EGT}. After the replacement, $\FF$ satisfies the $\Gamma(b)$-condition.

Set
\[
l_0 = \lf \frac{l \sqrt \ep}{\lambda} \rf, \hspace*{5mm}
\gamma = \frac{1}{\sqrt[4] \ep}, \eqand
b = \frac{4 \gamma m n}{l_0}.
\]
Then $\gamma$ is sufficiently large, and less than $\frac{l_0}{m^2}$ since
$
1 < \lambda < \frac{\ep l}{m^2}
$.

There exists a set $T$ such that $|T| \le \kappa \lp \FF \rp \Bigr/ \ln \frac{\ep l}{m^2 \lambda}$ and $\FF[T]$ satisfies the $\Gamma \lp b \rp$-condition in $X-T$: because the cardinality $j$ of such $T$ satisfies
\(
&&
|\FF| b^{-j} \le |\FF[T]| \le {n-j \choose m-j},
\\ &\Rightarrow&
b^{-j}\lp \frac{n}{m} \rp^j
\le
b^{-j} \prod_{j'=0}^{j-1} \frac{n-j'}{m-j'}
=
\frac{{n \choose m}}{{n-j \choose m-j}} b^{-j}
\le e^{\kappa \lp \FF \rp},
\\ &\Rightarrow&
j \le \frac{\kappa \lp \FF \rp}{\ln \frac{\ep l}{m^2 \lambda}}.
\)
Assume $j<m$, otherwise the desired claim is trivially true.

Apply \refco{EGT4Cor} to $\FF[T]$ noting
$\frac{l_0}{m^2} \le \frac{l_0-j}{(m-j)^2}$ and
$b \ge \frac{3 \gamma (m-j)(n-j)}{l_0-j}$.
We see 
\[
\left| Ext \lp \FF [T], l_0 \rp \right|
>
{n-j \choose l_0-j} \lp 1- \frac{2}{\sqrt[3] \gamma} \rp,
\]
from which
\[
\left| Ext \lp \FF[T], l \rp \right|
>
{n-j \choose l-j} \lp 1- e^{-\lambda} \rp,
\]
proving \refth{EGT}.  The truth of the last inequality is due to \reflm{PhaseII}: as $\frac{l - j}{l_0 -j} \ge \lambda \ep^{-1/2}$, it means
\[
\kappa \lbr {X \choose l} - Ext \lp \FF[T], l \rp \rbr
\ge 2^{\lf \log_2 \frac{l-j}{l_0-j} \rf} ~
\kappa \lbr {X \choose l_0} - Ext \lp \FF[T], l_0 \rp \rbr
> \lambda,
\]
in the universal set $X - T$ leading to the inequality.
This completes the proof of \refth{EGT}.

\section{About Splits of $X$}
Below we prove Theorems \ref{thSplit}, \ref{thSplit2} and \ref{thSplit3} given in Section 2.4.

\medskip 

{\noindent\bf \refth{thSplit}.} 
{\em Let 
\bdash 
\item a universal set $X$ of cardinality $n$ be weighted to induce a primitive norm $\| \cdot \|$, 
\item $m \in [n]$, and $r \in [m]$, 
\item $d_i \in [n]$ for $i \in [r]$ with $\sum_{i=1}^r d_i= n$, and $q_i \in [d_i]$ with $\sum_{i=1}^r q_i= m$. 
\edash 
For a family $\FF \subset {X \choose m}$, there exists a split $\bs{X}=(X_1, X_2, \ldots, X_r)$ of $X$ and subfamily $\FF' \subset \FF$ such that  
\[
\|\FF' \| \ge \frac{\|\FF \|}{{n \choose m}} \prod_{i=1}^r {d_i \choose q_i}, 
\quad 
|X_i|= d_i, 
\eqand 
|U \cap X_i| = q_i, 
\]
for all $i \in [r]$ and $U \in \FF'$. 
}
\begin{proof}
For such an $\FF \subset {X \choose m}$, $d_i$ and $q_i$, denote 
\[
\XX_j := \lb (X_1, X_2, \ldots, X_j)~:~  \textrm{$X_i \in {X \choose d_i}$ are mutually disjoint} \rb, 
\]
for $j \in [r]$, and 
\[
\FF_{\boldsymbol X} = \lb U~:~ U \in \FF, \textrm{~and~}
|U \cap X_i| = q_i \textrm{~for every~}  i \in [j] \rb, 
\] 
for $\bs{X} =(X_1, X_2, \ldots, X_j) \in \XX_j$. In addition, let $\XX_0 = \lb \emptyset \rb$ and $\FF_{\emptyset} = \FF$ for $j=0$. 

Our task here is to prove 
\beeq{eqthSplit}
\sum_{{\boldsymbol X} \in \XX_j} \left\| \FF_{\boldsymbol X} \right\| 
=
{n - \sum_{i=1}^j d_i   \choose m- \sum_{i=1}^j q_i} e^{-\kappa \lp \FF \rp}
\prod_{i=1}^{j}
{d_i \choose q_i} {n - \sum_{i'=1}^{i-1} d_{i'}  \choose d_i},
\eeq
for every $j \in [0, r-1]$. For if it is true for $j=r-1$, it means 
\(
\sum_{{\boldsymbol X} \in \XX_r} \left\| \FF_{\boldsymbol X} \right\| 
&=& 
\sum_{{\boldsymbol X} \in \XX_{r-1}} \left\| \FF_{\boldsymbol X} \right\| 
=e^{-\kappa \lp \FF \rp} \prod_{i=1}^{r} {d_i \choose q_i} {n - \sum_{i'=1}^{i-1} d_{i'}  \choose d_i}
\nexteqline 
\frac{\|\FF\|}{{n \choose m}} \prod_{i=1}^{r} {d_i \choose q_i} {n - \sum_{i'=1}^{i-1} d_{i'}  \choose d_i}. 
\)
As $|\XX_r| = \prod_{i=1}^{r} {n - \sum_{i'=1}^{i-1} d_{i'}  \choose d_i}$, it implies the existence of $\bs{X} \in \XX_r$ such that 
\\$\| \FF_{\bs{X}} \| \ge \frac{\|\FF\|}{{n \choose m}} \prod_{i=1}^{r} {d_i \choose q_i}$, proving the theorem.  
So \refeq{eqthSplit} suffices if confirmed for all $j$.

Proof of \refeq{eqthSplit} by induction on $j$ with the basis $j=0$ clearly true since 
\\$\sum_{{\boldsymbol X} \in \XX_0} \left\| \FF_{\boldsymbol X} \right\| 
=\| \FF \| = {n \choose m} e^{-\kappa(\FF)}$.

Assume true for $j$ and prove true for $j+1$. Fix any ${\boldsymbol X}  \in \XX_j$ putting
\(
&& 
X' = X - \bigcup_{i=1}^j X_i, \quad 
n' = |X'|, \quad 
\\ && 
m'= m-\sum_{i=1}^j q_i, \eqand
\gamma_{\boldsymbol X} = 
\frac{\left\| \FF_{\boldsymbol X} \right \|}
{{n' \choose m'} \prod_{i=1}^j {d_i \choose q_i} }. 
\)
Also write $d=d_{j+1}$, $q=q_{j+1}$ for simplicity, and 
$
{\boldsymbol X}' = \lp X_1, X_2, \ldots, X_j, X_{j+1} \rp 
$ 
when $X_{j+1} \in {X' \choose d}$ is given.

Consider pairs $\lp U, X_{j+1}\rp$ of $U \in \FF_{\boldsymbol X}$ and $X_{j+1}$ such that $\left| U \cap X_{j+1} \right|=q$. There are
${m' \choose q} {n'-m' \choose d-q}$ such pairs incident to each $U$. So the sum of $\| U  \|$ for all $\lp U, X_{j+1}\rp$ is 
\(
{m' \choose q}{n' - m' \choose d-q} \left\| \FF_{\boldsymbol X} \right\|
&=&
{m' \choose q} {n' - m' \choose d-q}  \gamma_{\boldsymbol X}   {n' \choose m'}  \prod_{i=1}^j {d_i \choose q_i}
\nexteqline
{n' \choose q} {n'-q \choose m'-q} {n' - m' \choose d-q}  \gamma_{\boldsymbol X}  \prod_{i=1}^j {d_i \choose q_i}. 
\)
Here 
\[
{n'-q \choose m'-q} {n' - m' \choose d-q} 
= \frac{(n'-q)!}{(m'-q)!(d-q)! (n-m'-d+q)!} 
= {n'-q \choose d-q} {n'-d \choose m'-q}. 
\]
So the above equals 
\[
{n' \choose q} {n'-q \choose d-q} {n'-d \choose m'-q} 
\gamma_{\boldsymbol X}  \prod_{i=1}^j {d_i \choose q_i} 
=
\gamma_{\boldsymbol X} 
{n'-d \choose m'-q}  {n' \choose d}
\prod_{i=1}^{j+1} {d_i \choose q_i}. 
\]

This finds the sum of $\| \FF_{\bs{X}'} \|$ for $\bs{X}'$ including the fixed $\bs{X}$. Now consider all $\bs{X} \in \XX_j$ to see 
\[
\sum_{{\boldsymbol X}' \in \XX_{j+1}} \left\| \FF_{\bs{X}'} \right\| 
= \sum_{{\boldsymbol X} \in \XX_j}  
\gamma_{\boldsymbol X} 
{n'-d \choose m'-q}  {n' \choose d}
\prod_{i=1}^{j+1} {d_i \choose q_i}. 
\] 
By induction hypothesis, $i.e.$ \refeq{eqthSplit} for $j$, we further notice 
\(
&&
\sum_{{\boldsymbol X} \in \XX_j} \left\| \FF_{\boldsymbol X} \right\|
=
{n' \choose m'} e^{-\kappa \lp \FF \rp}
\prod_{i=1}^{j}
{d_i \choose q_i} {n - \sum_{i'=1}^{i-1} d_{i'}  \choose d_i},
\sothat
\sum_{{\boldsymbol X} \in \XX_j}  \gamma_{\boldsymbol X}
=
e^{-\kappa \lp \FF \rp}
\prod_{i=0}^{j}
{n - \sum_{i'=1}^{i-1} d_{i'}  \choose d_i}. 
\)

By the above two, we conclude 
\[
\sum_{{\boldsymbol X}' \in \XX_{j+1}} \left\| \FF_{\bs{X}'} \right\| 
=
{n - \sum_{i=1}^{j+1} d_i   \choose m- \sum_{i=1}^{j+1} q_i} e^{-\kappa \lp \FF \rp}
\prod_{i=1}^{j+1}
{d_i \choose q_i} {n - \sum_{i'=1}^{i-1} d_{i'}  \choose d_i},
\]
confirming the induction step. We have proven \refeq{eqthSplit} for all $j$, hence the  theorem as well. 
\end{proof}

\medskip

{\noindent\bf \refth{thSplit2}.}  
{\em Let $X$ be weighted to induce a primitive norm $\| \cdot \|$, 
$m \in [n]$, and $r \in [m]$ such that $r \big| n$ and $\max \lp \frac{m}{n}, \frac{r}{m}, \frac{1}{r} \rp< \ep$. For a family $\FF \subset {X \choose m}$ with $\| \FF \| >0$, there exists an $r$-split $\bs{X}=(X_1, X_2, \ldots, X_r)$ and subfamily $\FF' \subset \FF$ satisfying 
\[
\|\FF' \| > \lbr 1-  r m \exp \lp - \frac{m}{12r} \rp \rbr \| \FF \|, \eqand \frac{m}{2r} < |U \cap X_i| < \frac{2m}{r}, 
\]
for all $U \in \FF'$ and strips $X_i$ of $\bs{X}$. 
}
\begin{proof} 
To prove the theorem, we first show that 
\beeq{eq0Split3}
\ln \frac{{m \choose q}{n-m \choose d - q}}{{n \choose d}}  < - \frac{m}{12 r}, 
\quad \forall q \in \left[ 0, \frac{m}{2r} \right] \cup \left[\frac{m}{2r}, m \right], 
\eqwhere d= \frac{n}{r}. 
\eeq 
For if $q=0$, 
\(
\frac{{n-m \choose d}}{{n \choose d}} &=& \prod_{i=0}^{d - 1} \lp 1- \frac{m}{n-i} \rp
< \lp 1- \frac{m}{n} \rp^{d} 
\\ &=& 
\lbr \lp 1- \frac{m}{n} \rp^{\frac{n}{m}}  \rbr^{\frac{m}{r}} 
<  \exp \lp - \frac{m}{2 r} \rp. 
\) 
If $q= d \le m$, we can similarly show $\ln {m \choose d} \Bigr/ {n \choose d} < - \frac{m}{r}$.

Else, confirm \refeq{eq0Split3} as follows given $q \not \in \lb 0, d \rb$. Use \reflm{SpaceProp} in Appendix C by setting 
\(
&& 
x=n, \quad x'=m, \quad y=\frac{n}{r}, \eqand y' = q, 
\sothat 
|\Delta| = \left| y' - \frac{x' y}{x} \right| \ge \frac{m}{2r}.
\) 
Also, 
\(
&& 
\beta = y' \lp 1- \frac{y'}{y} \rp \lp 1 - \frac{y-y'}{x-x'}\rp \lp 1- \frac{y'}{x} \rp \lp 1- \frac{y}{x} \rp^{-1}
> \lp 1- 2 \ep \rp^3, 
\sothat 
- \frac{\ln 2 \pi \beta}{2} + \frac{1}{12y} + \frac{1}{12(x-y)} < 0, 
\)
if $y'  \in [y-1]$: because $x ( 1- x/y)$ for the $y$ and any $x \in [y-1] $ takes the mininum $1- y^{-1}$ at $x=1$ and $y-1$. 

By the lemma and the two facts, we see 
$ 
\ln \frac{{m \choose q} {n-m \choose n/r -q}}{{n \choose n/r}}  < -\alpha \Delta \le - \frac{m}{12 r} 
$ confirming \refeq{eq0Split3}.

\medskip 

We now claim that for each $U$ and $i \in [r]$, there are less than $m n_s \exp \lp - \frac{m}{12r} \rp$ $r$-splits $\bs{X}$ such that 
\beeq{eqSplit3}
|U \cap X_i| \not \in \lp \frac{m}{2r},~ \frac{2m}{r} \rp, 
\eeq 
where $n_s$ is the total number of $r$-splits $\bs{X}$ given by 
\[
n_s := \prod_{i=0}^{r-1} {n - i d \choose d}. 
\]
This is true because for each $q \in [0, m]$, there are ${m \choose q}{n-m \choose d-q}$ sets $X_i \in {X \choose d}$ such that $|U \cap X_i| = q$. Then the claim holds because 
\[
{m \choose q}{n-m \choose d-q} < \exp \lp - \frac{m}{12r} \rp {n \choose d}, \quad \textrm{if}~  
q \not\in \lp \frac{m}{2r},~ \frac{2m}{r} \rp.  
\]
by \refeq{eq0Split3}.

Therefore, the sum of $\| U \|$ for all pairs $(U, \bs{X})$ such that \refeq{eqSplit3} for any $i$  is less than $rm n_s \exp \lp - \frac{m}{12r} \rp \| \FF \|$. By this we conclude the existence of an $\bs{X}$ and $\FF' \subset \FF$ meeting the desired two conditions, completing the proof. 
\end{proof}

\medskip 

{\noindent\bf \refth{thSplit3}.}
{\em Let $m \in [n]$, $r \in [m]$ such that 
\[
r \big| n, \quad 
\max \lp \frac{m}{n}, \frac{r}{m}, \frac{1}{r} \rp< \ep, \eqand b \in  \lp c^3 kr \ln \frac{m}{r}, \infty \rp, 
\]
for some constants $c \in ( \ep^{-1}, \infty)$ and $k \in \Z_{>0}$. For each nonempty $\FF \subset {X \choose m}$ satisfying the $\Gamma(b)$-condition, there exist an $r$-split $\bs{X}=(X_1, X_2, \ldots, X_r)$, $k$ mutually disjoint sets $C_i$, and subfamilies $\FF_i \subset \FF[C_i]  \cap {X - \bigcup_{i' \in [k] - \lb i \rb} C_{i'} \choose m}$ with the following four conditions satisfied for every $i \in [k]$: 
\be{i)}
\item $|C_i| < c r^2 \ln \frac{m}{r}$. 
\vspace*{1.2mm}
\item $|\FF_i|> \frac{1}{2} \lp \frac{r}{3m} \rp^r  b^{-|C_i|} |\FF|$. 
\item The $\Gamma\lp b / 6 \rp$-condition of $\FF_i$ in $X - C_i$. 
\item $|U \cap X_j| = q_{i, j}$ for every strip $X_j$ of $\bs{X}$, some $q_{i, j} \in \lp \frac{m}{2r},~ \frac{2m}{r}  \rp$, and all $U \in \FF_i$. 
\ee}
\begin{proof}
By \refth{thSplit2}, there exists an $r$-split $\bs{X}$, and subfamily $\FF'$ as described there with the unit norm $\| \cdot \| = |\cdot |$.  With them consider the following statement for $j \in [0, r]$.

\medskip 

{\em Proposition $\Psi(j)$:} there exist $k$ mutually disjoint sets $C_i$, and subfamilies $\FF_i \subset \FF'[C_i] \cap {X - \bigcup_{i' \in [k] - \lb i \rb} C_{i'} \choose m}$ satisfying the four conditions for every $i \in [k]$: 
\be{\hspace*{2mm} a)}
\item $|C_i| \le c j r \ln \frac{m}{r}$, where the inequality is strict if $j>0$. 
\vspace*{1.2mm}
\item $|\FF_i|> \frac{1}{2} \lp \frac{r}{3 m} \rp^j b^{-|C_i|} |\FF|$. 
\item $\Gamma\lp b_j, 2 \rp$-condition of $\FF_i$ in $X - C_i$, where $b_j = b \lp 1- \frac{1}{r} \rp^j$. 
\item $|U \cap X_{j'}| = q_{i, j'}$ for every strip $X_{j'}$ of $\bs{X}$ with $j' \in [j]$, some $q_{i, j'} \in \lp \frac{m}{2r},~ \frac{2m}{r}  \rp$, and all $U \in \FF_i$. 
\ee

Our remaining task is to prove $\Psi(j)$ on induction $j$ as $\Psi(r)$ implies the thoerem. For its basis $j=0$, choose  $\FF_i = \FF'$ and $C_i=\emptyset$ for all $i \in [k]$. It is straightforward\footnote{
Note $| \FF_i |>b^m$ since $1=|\FF[U] |<b^{-m} |\FF |$ for any $U \in \FF$ by the $\Gamma(b)$-condition of $\FF$. 
} to check that $\FF_i$, $C_i$ and $\bs{X}$ satisfy the four conditions $\Psi(0)$-a) to d). 
Assume $\Psi(j-1)$ for the induction step and prove $\Psi(j)$.

Let $i=1$, and $\FF_{i, q}$ for each $q \in [m]$ be the family of $U \in \FF_i$ such that $|U \cap X_j| = q$. Set $\FF'_i$ to be an $\FF_{i, q}$ with the largest cardinality over all $q$. 
As the $\bs{X}$ comes from \refth{thSplit2}, we have $q \in \lp \frac{m}{2r}, \frac{2m}{r}\rp$ so 
$
|\FF'_i| > \frac{r}{2m} |\FF_i| 
$.

Then find a maximal set $S_i \subset X - C_i$ such that 
\beeq{eq1Split3}
|\FF'_i[S_i]| \ge b_j^{-|S_i|} |\FF'_i|. 
\eeq 
Note the family $\FF'_i [S_i]$ satisfies the $\Gamma(b_j)$-condition in $X - C_i - S_i$.

In addition, 
\beeq{eq2Split3}
|S_i| < c r \ln \frac{m}{r}, 
\eeq
which is due to the $\Gamma(b_{j-1}, 2)$-condition of  $\FF_i$ from $\Psi(j-1)$-iii): it implies 
\(
&& 
\frac{r}{2m}  b_j^{-|S_i|} |\FF_i|
< 
|\FF'_i[S_i]| \le |\FF_i[S_i]|< 2 b_{j-1}^{-|S_i|}| \FF_i|,
\sothat 
\lp \frac{b_{j-1}}{b_j} \rp^{|S_i|} = \lp 1 -  \frac{1}{r} \rp^{-|S_i|} < \frac{4 m}{r}, 
\sothat 
|S_i| < c r \ln \frac{m}{r}. 
\)
Because if $|S_i| \ge c r \ln \frac{m}{r}$, we would have 
\(
\lp 1 -  \frac{1}{r} \rp^{-|S_i|} &\ge& \lp 1 -  \frac{1}{r} \rp^{-c r \ln \frac{m}{r}} 
= \lbr  \lp 1 -  \frac{1}{r} \rp^{-r} \rbr^{c \ln \frac{m}{r}} 
\\ &>& 
\lp \frac{e}{2} \rp^{c \ln \frac{m}{r}}> \frac{4 m}{r}, 
\)
contradicting $\lp 1 -  \frac{1}{r} \rp^{-|S_i|} < \frac{4 m}{r}$ above. This proves \refeq{eq2Split3}.

\medskip 

We now update $C_i$ and $\FF_i$ by 
\beeq{eq3Split3}
C_i \leftarrow C_i \cup S_i, \eqand  \FF_i \leftarrow \FF'_i[S_i], 
\eeq
for $i=1$. 
By $\Psi(j-1)$, \refeq{eq1Split3}, and \refeq{eq2Split3}, they satisfy $\Psi(j)$-a) and d), 
\beeq{eq3Split4}
|\FF_i|> \frac{3}{2} \lp \frac{r}{3 m} \rp^j b_j^{-|C_i|} |\FF|, 
\eeq 
instead of $\Psi(j)$-b), and the $\Gamma(b_j)$-condition in $X - C_i$ instead of $\Psi(j)$-c). 

\medskip 

Repeat the process for $i=2, 3, \ldots, k$ to update the other $C_i$ and $\FF_i$, noting the following: 
\bdash
\item Right after $C_1$ is updated by \refeq{eq3Split3}, we eliminate $\FF_i(C_1)$ from $\FF_i$ for each $i \in [2, k]$. This changes the $\Gamma\lp b_{j-1} \rp$-condition of $\FF_i$ in $X - C_i$ into the $\Gamma\lp b_{j-1},  1+ \frac{1}{ck} \rp$-condition: because $b_{j-1} > \ep c^3 k r \ln \frac{m}{r}$ and $|S_1| < c r \ln \frac{m}{r}$, the removal of $\FF_i(C_1)$ reduces $\FF_i$ by its less than $\frac{1}{c^2k}$.

\item In addition, right after $C_i$ is updated by \refeq{eq3Split3}, we eliminate $\FF_1(C_i)$ from the current $\FF_1$. This reduces less than $\frac{1}{c^2k}$ of $\FF_1$ as well.

\item We perform similar eliminations for each pair of $\FF_i$. The resulting recursive invariant is that $\FF_i$ satisfies the $\Gamma\lp b_{j-1},  1+ \frac{i-1}{ck} \rp$-condition in $X - C_i$ right before the update of $C_i$, and $\Gamma\lp b_j,  1+ \frac{k - i}{ck} \rp$-condition in $X - C_i$ right after the updates are over for all $C_{i'} \lp i' \in [k] \rp$. It is straightforward to check its truth arguing similarly to the above for $i=1$.

\item The process guarantees $\FF_i \subset \FF'[C_i] \cap {X - \bigcup_{i' \in [k] - \lb i \rb} C_{i'} \choose m}$, the mutual disjointness of $C_i$, and $\Psi(j)$-c) after its termination. 
\edash 

\medskip

With \refeq{eq3Split4}, those construct the final $\FF_i \subset \FF[C_i]$ with the four conditions of $\Pi(j)$ confirmed. We have proven the induction step to verify the theorem. 
\end{proof}

\section{Asymptotics of Binomial Coefficients} 
Stirling's approximation in form of a double inequality is known to be 
\beeq{Stirling}
\sqrt{2 \pi n} \lp \frac{n}{e} \rp^n
\exp \lp \frac{1}{12n+1} \rp
< n! <
\sqrt{2 \pi n} \lp \frac{n}{e} \rp^n
\exp \lp \frac{1}{12n} \rp,
\eeq
for any $n \in \Z_{>0}$ \cite{R55}. It is straightforward to derive from \refeq{Stirling} the following inequality on a binomial coefficient, called its {\em standard estimate}, $i.e.$,  
\beeq{eq1App}
\lp \frac{x}{y} \rp^y \le {x \choose y} < \lp \frac{e x}{y} \rp^y, 
\eeq 
for every $x, y \in \Z_{>0}$ with $x \ge y$.

We can further improve its accuracy. Let
\beeq{FunctionS}
s: (0, 1) \rightarrow (0,1),
\hspace*{5mm}
t \mapsto 1 - \lp 1 - \frac{1}{t} \rp \ln (1-t).
\eeq 
By the Taylor series of $\ln(1-t)$, the function is also expressed as
\beeq{ExpandS}
s(t)
=
1 +\lp 1 - \frac{1}{t} \rp \sum_{j \ge 1} \frac{t^j}{j}
= 1 + \sum_{j \ge 1} \frac{t^j}{j} - \sum_{j\ge 0} \frac{t^j}{j+1}
= \sum_{j \ge 1} \frac{t^j}{j(j+1)}.
\eeq
With this we can show the following double inequality as in \cite{3SF}: 

\begin{lemma} \label{asymptotic} 
For $x, y \in \Z_{>0}$ such that $x<y$,
\(
&&
\frac{1}{12 x+1} - \frac{1}{12 y} - \frac{1}{12(x-y)}
<
z
<
\frac{1}{12 x} - \frac{1}{12 y+1} - \frac{1}{12(x-y)+1},
\\ \textrm{where} &&
z = \ln {x \choose y} - y \lbr \ln \frac{x}{y} + 1 - s \lp \frac{y}{x} \rp  \rbr
- \frac{1}{2} \ln \frac{x}{2 \pi y(x-y)}.  
\)
\end{lemma}
\begin{proof}
By \refeq{Stirling} find that 
\(
&&
\frac{1}{12 x+1} - \frac{1}{12 y} - \frac{1}{12(x-y)}
<
\ln \frac{{x \choose y}}{u}
<
\frac{1}{12 x} - \frac{1}{12 y+1} - \frac{1}{12(x-y)+1},
\\&&
\eqwhere
u
= \sqrt{\frac{x}{2 \pi y(x-y)}} \frac{x^x}{y^y (x-y)^{x-y}}.
\)
Since
$
\ln \frac{x^x}{y^y (x-y)^{x-y}}  = y \lbr \ln \frac{x}{y} + 1 - s \lp \frac{y}{x} \rp  \rbr
$
by \refeq{FunctionS}, it proves the lemma.
\end{proof}

Hence, 
\beeq{eqAsymptotic} 
\ln {x \choose y}  = y \lbr \ln \frac{x}{y} + 1 - s \lp \frac{y}{x} \rp  \rbr
+  \frac{1}{2} \ln \frac{x}{2 \pi y(x-y)} + z, 
\eeq 
for any $x, y \in \Z_{>0}$ with $y<x$, and some 
\[
z\in \lp - \frac{1}{12y} - \frac{1}{12(x-y)} ,~0 \rp. 
\]
Confirm the following lemma with \refeq{eqAsymptotic} that is referred to by our proof of \refth{thSplit2} in Appendix B. 

\begin{lemma} \label{SpaceProp} 
For $x \in \Z_{>1}$, $y \in [x-1]$, $x' \in [1, x/2] \cap \Z$, and $y' \in [x'-1] \cap [y-1]$, 
\(
&&
\ln \frac{{x-x' \choose y-y'}{x' \choose y'}}{{x \choose y} } 
< 
- \alpha \Delta
-  \frac{\ln 2 \pi \beta}{2} + \frac{1}{12y} +  \frac{1}{12(x-y)}, 
\\ \textrm{where} && 
\Delta = y' - \frac{x' y}{x}, \quad 
\alpha = \lb \begin{array}{cc} 
x \Delta \big/ (3x'y)  & \textrm{if $y' / x' < 2 y / x$,} \\
\ln 4 - 1   & \textrm{otherwise,} \\
\end{array} \right. 
\\ \textrm{and} && 
\beta = y' \lp 1- \frac{y'}{y} \rp \lp 1- \frac{y- y'}{x- x'} \rp \lp 1- \frac{y'}{x'} \rp 
\lp 1 - \frac{y}{x} \rp^{-1}.  
\) 
\end{lemma}
\begin{proof}
Evaluate the three binomial coefficients by \refeq{eqAsymptotic} to find  
\[
\ln {x-x' \choose y-y'}{x' \choose y'} - \ln {x \choose y}
< u - \sum_{j \ge 1} \frac{v_j}{j(j+1)} 
- \frac{\ln 2 \pi \beta}{2} + \frac{1}{12y} + \frac{1}{12(x-y)}, 
\]
where
\(
&&
u =
\lp y-  y' \rp \lp \ln \frac{x-x'}{y-y'} + 1\rp
+ y' \lp \ln \frac{x'}{y'} + 1 \rp
- y \lp \ln \frac{x}{y} + 1 \rp,
\\ \textrm{and} &&
v_j =\frac{(y-y')^{j+1}}{(x-x')^j} + \frac{y'^{~j+1}}{x'^{~j}}
- \frac{y^{j+1}}{x^j}. 
\) 
So it suffices to show 
\beeq{eq2SpaceProp} 
v_j \ge 0, \quad \textrm{for all~} j \in \Z_{>0}, 
\eeq
and 
\beeq{eqSpaceProp}
u< -\alpha \Delta, \quad \textrm{if}~ \Delta\ne 0. 
\eeq

To see \refeq{eq2SpaceProp}, put 
\[
a = \frac{x'}{x}, \quad 
b =\frac{y'}{y}, \eqand 
f =
\frac{\lp 1- b \rp^{j+1}}{\lp 1- a \rp^j}
+ \frac{b^{j+1}}{a^j}, 
\]
for a given $j$. 
Then the desired condition holds if $f \ge 1$ for each fixed $a \in (0, 1)$ and all $b \in [0, 1]$. It is straightforward to check its truth with basic calculus.

We show \refeq{eqSpaceProp} by finding that 
\(
u &=&
\lp y-  y' \rp \lp \ln \frac{x}{y} + \ln \frac{1-\frac{x'}{x}}{1-\frac{y'}{y}} \rp
+ y' \ln \frac{x'}{y'} 
- y  \ln \frac{x}{y} 
\\ &=& 
y' \ln \frac{x'y}{xy'} + \lp y - y' \rp \ln \frac{1- \frac{x'}{x}}{1-\frac{y'}{y}}.  
\)
Then put 
\[
p = \frac{x \Delta}{x'y}, \eqand
q= \frac{\Delta}{y \lp 1-\frac{x'}{x} \rp}. 
\]
to derive 
\(
&& 
\frac{x'y}{x}  = \frac{\Delta}{p}, \quad 
y' = \Delta \lp \frac{1}{p} +1 \rp, \quad 
1-q=\frac{1-\frac{y'}{y}}{1- \frac{x'}{x}}, \quad 
\Delta \lp \frac{1}{q} - 1 \rp = y-y', 
\sothat 
u=
- \Delta \lp \frac{1}{p} + 1 \rp \ln \lp 1 + p \rp
- \Delta \lp \frac{1}{q} -1 \rp  \ln \lp 1 -  q \rp. 
\)

Observe some facts.
\bdash 
\item $p > -1$ otherwise $\frac{xy'}{x'y}  \le 0$. Also $p \ge 1$ $\Leftrightarrow$ $\frac{y'}{x'} \ge \frac{2 y}{x}$. 

\item $|q| < 1 \Leftrightarrow \frac{2x'}{x} -1 < \frac{y'}{y} < 1$ is true by $x' \le x/2$.

\item By the Taylor series of natural logarithm,
\(
&& \lp \frac{1}{p} + 1 \rp \ln \lp 1 + p \rp
= 1 - \sum_{j \ge 1} \frac{(-p)^j}{j(j+1)},
\quad \textrm{if $|p|<1$,}
\\ \textrm{and} &&
\lp \frac{1}{q} -1 \rp  \ln \lp 1 -  q \rp
= - 1 + \sum_{j \ge 1} \frac{q^j}{j(j+1)},
\quad \textrm{as $|q|<1$.}
\)

\item Note $p \Delta$ and $q \Delta$ are both positive, so 
\(
&&
-\Delta \lp \frac{1}{p} + 1 \rp \ln \lp 1 + p \rp
< \Delta \lp -1 - \frac{p}{2} + \frac{p^2}{6}  \rp, \quad \textrm{if $|p|<1$,}
\\ \textrm{and} && 
- \Delta \lp \frac{1}{q} -1 \rp  \ln \lp 1 -  q \rp < \Delta. 
\) 
The first inequality is due to $-\Delta \lbr\frac{p^j}{j(j+1)} - \frac{p^{j+1}}{(j+1)(j+2)} \rbr < 0$ for all odd $j \ge 3$, and the second inequality to $-\Delta \lbr \frac{q^j}{j(j+1)} - \frac{q^{j+1}}{(j+1)(j+2)} \rbr < 0$ for odd $j \ge 1$. 
\edash

\medskip

Therefore, 
\[
u < \Delta \lp \frac{-p}{2} + \frac{p^2}{6} \rp < - \frac{p \Delta}{3} = - \alpha \Delta, 
\]
if $|p|<1$. We now have \refeq{eqSpaceProp} when $\frac{y'}{x'} < \frac{2 y}{x}$.

If $p \ge 1$, then $\lp 1 + \frac{1}{p} \rp \ln \lp 1 + p \rp \ge 2 \ln 2$, since the LHS monotonically increases in $p$. Hence, 
\[
u < \lp -2 \ln 2 + 1  \rp \Delta = - \alpha \Delta, 
\]
in this case as well. This completes the proof. 
\end{proof}

\bibliographystyle{amsplain}

\end{document}